\documentclass[smallextended,envcountsect]{svjour3}

\usepackage{amsmath, amssymb, dsfont, float, graphicx, anysize}

\newcommand\I{\mathcal{I}}
\newcommand\charfun{\mathds{1}}
\renewcommand\le{\leqslant}
\renewcommand\ge{\geqslant}
\newcommand{\EndResult}{\hspace*{0pt}\hfill $\square$}

\title{Intriguing sets of partial quadrangles}

\author{John Bamberg \and Frank De Clerck \and Nicola Durante}
\institute{ %
John Bamberg \at
Department of Pure Mathematics and Computer Algebra,
Ghent University,
Krijgslaan 281--S22, B-9000 Ghent, Belgium.
\email{bamberg@cage.ugent.be}
\and
Frank De Clerck \at
Department of Pure Mathematics and Computer Algebra,
Ghent University,
Krijgslaan 281--S22, B-9000 Ghent, Belgium.
\email{fdc@cage.ugent.be}
\and
Nicola Durante
\at
Dipartimento di Matematica ed Applicazioni,
Universit\`a di Napoli ``Federico II'',
80125 Naples,
Italy.
\email{ndurante@unina.it}
}

\begin{document}

\maketitle


\begin{abstract}
  The point-line geometry known as a \textit{partial quadrangle} (introduced by
  Cameron in 1975) has the property that for every point/line non-incident pair
  $(P,\ell)$, there is at most one line through $P$ concurrent with $\ell$. So
  in particular, the well-studied objects known as \textit{generalised
    quadrangles} are each partial quadrangles.  An \textit{intriguing set} of a
  generalised quadrangle is a set of points which induces an equitable partition
  of size two of the underlying strongly regular graph. We extend the theory of
  intriguing sets of generalised quadrangles by Bamberg, Law and Penttila to
  partial quadrangles, which surprisingly gives insight into the structure of
  hemisystems and other intriguing sets of generalised quadrangles.
 \end{abstract}


\section{Introduction}

A set of points $\I$ of a generalised quadrangle is defined in \cite{BLP} to be
\textit{intriguing} if the number of points of $\I$ collinear to an arbitrary
point $P$ is a constant $h_1$ if $P$ lies in $\I$, and another constant $h_2$ if
$P$ resides outside of $\I$. For example, a line of a generalised quadrangle is
such an object where $h_1$ is the number of points on a line, and
$h_2=1$. Eisfeld \cite{Eisfeld99} asks whether such sets have a natural
geometric interpretation, and it is shown in \cite{BLP} that the intriguing sets
of a generalised quadrangle are precisely the $m$-ovoids and tight sets
introduced by J. A. Thas \cite{Thas89} and S. E. Payne \cite{Payne87}
respectively. If one looks to the point graph of a generalised quadrangle, one
will find a strongly regular graph.  The associated Bose-Mesner algebra of this
graph decomposes into an orthogonal decomposition of three eigenspaces of the
adjacency matrix, one of which is the one-dimensional subspace generated by the
``all $1$'s'' vector. The other two eigenspaces correspond naturally to the two
types of intriguing sets; the positive eigenvalue corresponds to the tight sets,
and the negative eigenvalue corresponds to the $m$-ovoids \cite[Theorem
4.1]{BLP}.  In the broader context of association schemes, the
\textit{intriguing sets} correspond to the $\{0,1\}$ valued elements of the
Bose-Mesner algebra which are annihilated by all but one of the nontrivial
minimal idempotents (n.b., the trivial minimal idempotent has rank $1$).  So one
can employ the same techniques and exploit the orthogonal decomposition of the
associated Bose-Mesner algebra to derive information about certain geometric
configurations (see \cite{BKLP}, \cite{DeWVM} and \cite{Eisfeld99}).  In this
paper we consider the algebraic combinatorics of a partial quadrangle, a
geometric object which comes equipped with an interesting Bose-Mesner algebra.

A \textit{partial quadrangle} was introduced by P. J. Cameron \cite{Cameron75}
as a geometry of points and lines such that every two points are on at most one
line (and hence two lines meet in at most one point), there are $s+1$ points on
a line, every point is on $t+1$ lines and satisfying the following two important
properties:
\begin{enumerate}
\item[(i)] for every point $P$ and every line $\ell$ not incident with $P$,
  there is at most one point on $\ell$ collinear with $P$;
\item[(ii)] there is a constant $\mu$ such that for every pair of non-collinear
  points $(X,Y)$ there are precisely $\mu$ points collinear with $X$ and $Y$.
\end{enumerate}
With the above specifications, we say that the partial quadrangle has
\emph{parameters} $(s,t,\mu)$, or that it is a partial quadrangle
$\mathsf{PQ}(s,t,\mu)$. Note that the point-graph of this object is strongly
regular (see Section \ref{alggraph}).

The only known partial quadrangles, which are not generalised quadrangles, are
\begin{itemize}
\item triangle-free strongly regular graphs (i.e., partial quadrangles with
  $s=1$);
\item one of three exceptional examples, namely they arise from linear
  representation of one of the Coxeter $11$-cap of $\mathsf{PG}(4,3)$, the Hill
  $56$-cap of $\mathsf{PG}(5,3)$ or the Hill $78$-cap of $\mathsf{PG}(5,4)$;
\item or arise from removing points from a generalised quadrangle of order $(s,s^2)$.
\end{itemize}
We will now be more precise for this last class of partial quadrangles.  Let
$\mathcal{G}$ be a generalised quadrangle of order $(s,s^2)$ and let $P$ be a
point of $\mathcal{G}$. Then by removing all those points $P^\perp$ which are
collinear with $P$ results in a partial quadrangle $\mathsf{PQ}(s-1,s^2,s(s-1))$
(see \cite[pp. 4]{CDG79}). We will often refer to this construction as a
generalised quadrangle minus the perp of a point. Similarly, we can remove a
certain type of $m$-ovoid from $\mathcal{G}$ to obtain a partial quadrangle
\cite[Prop. 2.2]{CDG79}. A \textit{hemisystem} of $\mathcal{G}$, where $s$ is
odd, is a set of points $\mathcal{H}$ of $\mathcal{G}$ such that every line
meets $\mathcal{H}$ in $(s+1)/2$ points (i.e., it is an $m$-ovoid with
$m=(s+1)/2$). By considering the incidence structure restricted to
$\mathcal{H}$, we obtain a partial quadrangle
$\mathsf{PQ}((s-1)/2,s^2,(s-1)^2/2)$. Recently, Cossidente and Penttila
\cite{CossidentePenttila05} have constructed new hemisystems of the classical
generalised quadrangle $\mathsf{Q}^-(5,q)$, and thus new partial quadrangles.
In \cite{BDCD}, a hemisystem is constructed of the dual of the
Fisher-Thas-Walker-Kantor generalised quadrangle of order $(5,5^2)$, yielding a
new $\mathsf{PQ}(2,25,8)$.

For generalised quadrangles, it has been shown that an $m$-ovoid and an
$i$-tight set intersect in $mi$ points \cite[Theorem 4.3]{BLP}. From this
observation, one can prove or reprove interesting results in the forum of
generalised quadrangles. For partial quadrangles, the theory still holds; there
are two types of intriguing sets according to the parity of the associated
eigenvalue, and there is a similar ``intersection result'' (see Section
\ref{general}).  In Section \ref{thin}, we investigate and in some cases
classify, the intriguing sets of triangle-free strongly regular graphs; the
\textit{thin} partial quadrangles.  The section that follows concerns the
two known families of \textit{thick} partial quadrangles which arise from (i)
deleting the perp of a point, or from (ii) deleting a hemisystem. In both cases,
we look to the deleted point set, which we nominate as ``infinity'', and analyse
the situation for when an intriguing set of the ambient generalised quadrangle
gives rise to an intriguing set of the partial quadrangle obtained by removing
infinity.  In the case of a generalised quadrangle minus the perp of a point, we
give some strong combinatorial information in Section \ref{GQminusperp} on the
structure of incumbent intriguing sets, which manifests in a characterisation of
the \emph{positive} intriguing sets arising from tight sets of the ambient
generalised quadrangle, and a partial characterisation of the \emph{negative}
intriguing sets. The intriguing sets of partial quadrangles obtained from
hemisystems have less combinatorial structure, however, we are able to deduce
certain relationships between intriguing sets of the ambient generalised
quadrangle and the partial quadrangle (see Section \ref{GQminushemi}).  In
Section \ref{linear}, we return to isolated examples of partial quadrangles, and
this time on the exceptional examples arising from caps of projective spaces via
linear representation.

We will next revise and introduce the necessary material from algebraic graph
theory, including strongly regular graphs, the Bose-Mesner algebra and minimal
idempotents.  The notion of an intriguing set of a partial quadrangle then
follows from the more natural setting of an intriguing set of a strongly regular
graph. This allows us to focus our concentration on the combinatorics of the
underlying graph and its Bose-Mesner algebra.


\section{Some algebraic graph theory and intriguing sets}\label{alggraph}

\subsection{Intriguing sets of strongly regular graphs}

A regular graph $\Gamma$, with $v$ vertices and valency $k$, is \emph{strongly
  regular} with parameters $(v,k,\lambda,\mu)$ if (i) any two adjacent vertices
are both adjacent to $\lambda$ common vertices; (ii) any two non-adjacent
vertices are both adjacent to $\mu$ common vertices. If $A$ is the adjacency
matrix of the strongly regular graph $\Gamma$, then $A$ has three eigenvalues
and satisfies the equation $A^2=kI+\lambda A+\mu(J-I-A)$ where $I$ is the
identity matrix and $J$ is the all-ones matrix.  The all-ones vector $\charfun$
is an eigenvector of $A$ with eigenvalue $k$. The remaining two eigenvalues
$e^+$ and $e^-$ satisfy the quadratic equation $x^2=k+\lambda
x+\mu(-1-x)$. Hence $\mu-k=e^+e^-$ and $\lambda-\mu=e^++e^-$. (Since $A$ has $0$
trace, we deliberately write $e^+$ and $e^-$ in accordance with their parity;
one is nonnegative, and the other is negative).

As mentioned in the introduction, a strongly regular graph comes equipped with
its Bose-Mesner algebra, the $3$-dimensional matrix algebra generated by $A$,
$I$ and $J$. Now the Bose-Mesner algebra of a strongly regular graph is a
commutative algebra of real symmetric matrices, and so it has an orthogonal
decomposition into idempotents.  By idempotent, we mean with respect to ordinary
matrix multiplication. Moreover, there exist so called \textit{minimal
  idempotents} $E_0, E_1, E_2$ such that the product of any two is zero, and
such that they add up to the identity matrix. To obtain these matrices, one can
take the Gram matrices of the orthogonal projections to the three eigenspaces of
$A$.  So for a strongly regular graph with eigenvalues $k$ (the valency), $e^+$
and $e^-$, we can take the following minimal idempotents (n.b., $n$ is the size
of $A$):
\begin{align*}
E_0&=\frac{1}{n} J,\\
E_1&=\frac{1}{e^+-e^-}\left(A - e^-I - \frac{k-e^-}{n}J\right),\\
E_2&=\frac{1}{e^--e^+}\left(A - e^+I - \frac{k-e^+}{n}J\right).\\
\end{align*}
All of the above content is standard in the theory of association schemes and
can be found in a textbook such as \cite{Godsil93}.

We say that a set of vertices $\mathcal{I}$ of a strongly regular graph $\Gamma$
is an \textit{intriguing set} with parameters $(h_1,h_2)$ if there are two
constants $h_1$ and $h_2$ such that the number of elements of $\mathcal{I}$
adjacent to any vertex of $\mathcal{I}$ is $h_1$, and the number of elements of
$\mathcal{I}$ adjacent to any vertex of $\Gamma\setminus\mathcal{I}$ is $h_2$.
So necessarily, the subgraph induced by $\mathcal{I}$ is regular of valency
$h_1$.  We will call $h_1$ and $h_2$ the \textit{intersection numbers} of $\I$,
and note that we have made a slight difference here in comparison to the
definition in \cite{BLP}; our parameter $h_1$ will always be one less than the
analogue in \cite{BLP} due to ``adjacency'' being an anti-reflexive relation. It
turns out (see Lemma \ref{eigen}) that $h_1-h_2$ is an eigenvalue of the
adjacency matrix, and so we define $\mathcal{I}$ to be a \textit{positive} or
\textit{negative} intriguing set according to whether $h_1-h_2$ is equal to
$e^+$ or $e^-$.

From the algebraic graph theoretic point of view, an intriguing set of a
strongly graph is a set of vertices whose characteristic vector is annihilated
by one of the minimal idempotents $E_1$ or $E_2$. This simple observation allows
us to design algorithms to search for intriguing sets. A characteristic vector
of a set of points has values $0$ or $1$, so an intriguing set corresponds to a
set of rows of a minimal idempotent which add to the zero vector. One can reduce
the problem by taking the row Echelon reduced form of the given minimal
idempotent or by using subgroups of the induced permutation group on the points
to obtain collapsed matrices with constant row sums.

The following results follow in the same way as in \cite{BKLP} (see also
\cite{Eisfeld99}). We will assume throughout this paper that the entire vertex
set is not an intriguing set, and hence, that $h_1\ne h_2$. We use the notation
$\charfun_\I$ for the characteristic vector of $\I$.

\begin{lemma} \label{eigen} 
  Let $\I$ be an intriguing set of a
  strongly regular graph $\Gamma$, and let the intersection numbers of $\I$ be
  $h_1$ and $h_2$. Let $v$ and $k$ be the number of vertices and the valency of
  $\Gamma$ respectively, and let $A$ be the adjacency matrix of $\Gamma$.  Then:
  \begin{enumerate}
  \item[(i)] $(h_1 - h_2 - k ) \charfun_\I + h_2 \charfun$ is an
  eigenvector of $A$ with eigenvalue $h_1 - h_2 $;
  \item[(ii)] $|\I|=h_2v/(k-h_1+h_2).$
  \end{enumerate}
\end{lemma}

\begin{proof}
  The proof of (i) is just a straight-forward calculation, so we provide the
  proof for part (ii).  Let $A$ be the adjacency matrix of $\Gamma$.  Since $A$
  is a real symmetric matrix, the eigenvector $(h_1 - h_2 - k) \charfun_\I + h_2
  \charfun$ is orthogonal to the all-ones vector $\charfun$ with eigenvalue
  $k$.  So $\left((h_1 - h_2 - k) \charfun_\I + h_2 \charfun\right)\cdot
  \charfun=0$ and
  hence: $$-(h_1-h_2-k)\charfun_\I\cdot\charfun=h_2\charfun\cdot\charfun$$ from
  which the conclusion follows.
\EndResult\end{proof}

\begin{lemma}[Intersection Lemma]\label{intersection}
  Let $\I^+$ and $\I^-$ be positive and negative intriguing sets respectively of
  a strongly regular graph $\Gamma$ and let $v$ be the total number of
  vertices. Then
$$|\I^+\cap\I^-|= |\I^+||\I^-|/v .$$
\end{lemma}

\begin{proof}
  Just as in \cite[Theorem 4]{BKLP}, we use the fact that the eigenvectors
  corresponding to $\I^+$ and $\I^-$ (see Lemma \ref{eigen}) are orthogonal from
  which the result easily follows.
  \EndResult\end{proof}

One can obtain new intriguing sets by taking unions of disjoint intriguing sets
of the same type. Moreover, the complement of an intriguing set is also
intriguing, and of the same type. These observations will be important in
the study of intriguing sets of strongly regular graphs.

\begin{lemma}\label{newfromold}
  Suppose we have a strongly regular graph $\Gamma$ and let $\mathcal{A}$ and
  $\mathcal{B}$ be two intriguing sets of the same type, that is, they give rise
  to eigenvectors with the same eigenvalue. Then:
\begin{enumerate}
\item[(a)] If $\mathcal{A}\subset \mathcal{B}$, then $\mathcal{B}\setminus
  \mathcal{A}$ is an intriguing set of the same type as $\mathcal{B}$ and
  $\mathcal{A}$;
\item[(b)] If $\mathcal{A}$ and $\mathcal{B}$ are disjoint, then
  $\mathcal{A}\cup \mathcal{B}$ is an intriguing set of the same type as
  $\mathcal{A}$ and $\mathcal{B}$; 
\item[(c)] The complement $\mathcal{A}'$ of $\mathcal{A}$ in $\Gamma$
   is an intriguing set of the same type as $\mathcal{A}$.
\end{enumerate}
\end{lemma}

\begin{proof}
  Let $\mathcal{A}$ and $\mathcal{B}$ have intersection numbers $(a_1,a_2)$ and
  $(b_1,b_2)$ respectively. For part (a), note that if $P$ is a vertex in
  $\mathcal{B}\setminus\mathcal{A}$, then there are $b_1-a_2$ vertices of
  $\mathcal{B}\setminus\mathcal{A}$ adjacent to $P$. On the other hand, if $P$
  is in the complement of $\mathcal{B}\setminus\mathcal{A}$, that is, either in
  $\mathcal{A}$ or in the complement of $\mathcal{B}$, then there are $b_1-a_1$
  or $b_2-a_2$ vertices of $\mathcal{B}\setminus\mathcal{A}$ adjacent to $P$
  accordingly. However, since $\mathcal{A}$ and $\mathcal{B}$ are of the same
  type, we have that $b_1-b_2=a_1-a_2$ and hence the quantities $b_2-a_2$ and
  $b_1-a_1$ are equal. Therefore $\mathcal{B}\setminus\mathcal{A}$ is
  intriguing.  For part (b), it is simple to calculate that there are
  $a_1+b_2=a_2+b_1$ vertices of $\mathcal{A}\cup\mathcal{B}$ adjacent to $P$ if
  $P\in\mathcal{A}\cup\mathcal{B}$, but $a_2+b_2$ when
  $P\notin\mathcal{A}\cup\mathcal{B}$. Thus $\mathcal{A}\cup\mathcal{B}$ is
  intriguing and of the same type as $\mathcal{A}$ and $\mathcal{B}$ (by knowing
  the difference of the intersection numbers).  For part (c), let $k$ be the
  valency of the strongly regular graph. Then clearly, there are $k-a_2$
  neighbours of $P$ in $\mathcal{A}'$ if $P\in\mathcal{A}'$, but $k-a_1$
  neighbours in $\mathcal{A}'$ when $P\in\mathcal{A}$. Thus the complement of
  $\mathcal{A}$ is intriguing and of the same type as $\mathcal{A}$.
\EndResult\end{proof}


%


Below we give a simple example of how to determine (by hand) the intriguing sets
of the Petersen graph.


\subsection{Example: Intriguing sets of the the Petersen graph}

The two minimal idempotents we will consider of the Petersen graph are:

$$E_1=\tfrac{1}{6}\left(\tiny{\begin{array}{rrrrrrrrrr}
     3&   1&  -1&  -1&   1&   1&  -1&  -1&  -1&  -1 \\
     1&   3&   1&  -1&  -1&  -1&  -1&   1&  -1&  -1 \\
    -1&   1&   3&   1&  -1&  -1&  -1&  -1&  -1&   1 \\
    -1&  -1&   1&   3&   1&  -1&   1&  -1&  -1&  -1 \\
     1&  -1&  -1&   1&   3&  -1&  -1&  -1&   1&  -1 \\
     1&  -1&  -1&  -1&  -1&   3&   1&  -1&  -1&   1 \\
    -1&  -1&  -1&   1&  -1&   1&   3&   1&  -1&  -1 \\
    -1&   1&  -1&  -1&  -1&  -1&   1&   3&   1&  -1 \\
    -1&  -1&  -1&  -1&   1&  -1&  -1&   1&   3&   1 \\
    -1&  -1&   1&  -1&  -1&   1&  -1&  -1&   1&   3 
  \end{array}}\right),\quad
E_2=\tfrac{1}{15}\left(\tiny{\begin{array}{rrrrrrrrrr}
   6&  -4&   1&   1&  -4&  -4&   1&   1&   1&   1 \\
    -4&   6&  -4&   1&   1&   1&   1&  -4&   1&   1 \\
     1&  -4&   6&  -4&   1&   1&   1&   1&   1&  -4 \\
     1&   1&  -4&   6&  -4&   1&  -4&   1&   1&   1 \\
    -4&   1&   1&  -4&   6&   1&   1&   1&  -4&   1 \\
    -4&   1&   1&   1&   1&   6&  -4&   1&   1&  -4 \\
     1&   1&   1&  -4&   1&  -4&   6&  -4&   1&   1 \\
     1&  -4&   1&   1&   1&   1&  -4&   6&  -4&   1 \\
     1&   1&   1&   1&  -4&   1&   1&  -4&   6&  -4 \\
     1&   1&  -4&   1&   1&  -4&   1&   1&  -4&   6 
\end{array}}\right).$$

To obtain the intriguing sets, we first look for rows of $E_1$ which add to the
zero vector. We will identify the vertices of the Petersen graph, and hence the
rows of $E_1$, with $\{1,2,\ldots,10\}$.  Since the Petersen graph is vertex
transitive, we may suppose without loss of generality that $1$ is in our
putative intriguing set $\I$.  It turns out that the stabiliser of $1$ in the
automorphism group of the Petersen graph has as orbits $\{2,5,6\}$ and
$\{3,4,7,8,9,10\}$. We can see this by looking at the values in the first column
of $E_1$. So we may suppose without loss of generality that $3\in\I$. So far,
our two rows of $\I$ add to $\tfrac{1}{6}(2, 2, 2, 0, 0, 0, -2, -2, -2, 0)$. It
turns out that $\I$ can only be one of $\{ 1, 3, 4, 6, 8, 9\}$, $\{1, 3, 5, 7,
8, 10\}$ or $\{1, 3, 7, 9\}$; and we can exclude the first two since the
collection of intriguing sets are closed under complements and hence we can
regard only those of size at most $5$. The set $\{1,3,7,9\}$ corresponds to a
4-coclique of the Petersen graph. In fact, there are in total five 4-cocliques
of the Petersen graph.

For the second minimal idempotent $E_2$, we similarly assume that $1$ and $2$
are contained in our putative intriguing set $\I$. The sum of the first two rows
of $E_2$ is $\tfrac{1}{15}(2, 2, -3, 2, -3, -3, 2, -3, 2, 2)$, and in order to
cancel this vector, we must complete $\I$ to one of $\{1,2, 3, 4, 5\}$, $\{1,2,
3, 6, 10\}$, $\{1,2,5, 8, 9\}$, or $\{1,2, 6, 7, 8\}$. It then follows that the
twelve pentagons (5-cycles) are intriguing sets of the Petersen graph.


\subsection{Intriguing sets of partial quadrangles, the basics}\label{general}

Let $\mathcal{P}$ be a point-line incidence structure whose point graph is
strongly regular. Then a set of points of $\mathcal{P}$ is an \textit{intriguing
  set} if it corresponds to an intriguing set of the point graph. We will use
the symbol $\perp$ to denote the collinearity relation on points, so $P^\perp$
will denote the set of all points collinear to $P$. However, we will also extend
the graph theoretic notion of adjacency to geometries by writing $P^\sim$ to
mean the set of all points collinear \textbf{but not equal to} $P$; that is, the
\textit{neighbours} of $P$.  (Thus our point graphs have no loops, and our
adjacency matrices will have $0$'s on the diagonal).

Let $\mathcal{G}$ be a generalised quadrangle of order $(s,t)$.  The point graph
of $\mathcal{G}$ is strongly regular with parameters:
$$v=(s+1)(st+1), \quad k=s(t+1),\quad \lambda=s-1,\quad\mu=t+1,$$
and hence has three eigenvalues, one of which is the valency $k$.  The other two
eigenvalues, commonly known as the \textit{principal eigenvalues}, are $s-1$ and
$-t-1$. The eigenvalues of the point graph of a partial quadrangle with
parameters $(s,t,\mu)$ are accordingly:

\begin{center}
\begin{tabular}{l|c}
the valency& $s(t+1)$\\
positive & $e^+:=(-\mu-1+s+\sqrt{(\mu-1-s)^2+4st})/2$\\
negative & $e^-:=(-\mu-1+s-\sqrt{(\mu-1-s)^2+4st})/2$.
\end{tabular}
\end{center}


From the above definition, a nonempty subset of points $\I$ of a partial
quadrangle $\mathsf{PQ}(s,t,\mu)$ is intriguing if there are two constants $h_1$
and $h_2$ such that
$$|P^\sim\cap\I|=\begin{cases}
h_1&\text{if }P\in \I,\\
h_2&\text{otherwise}\\
\end{cases}$$ where $P$ runs over the points of the partial quadrangle.  In
other words, if $A$ is the adjacency matrix of the point graph, then $\I$ is
intriguing if and only if its characteristic function $\charfun_\I$ satisfies the
following relation:
$$A\charfun_\I=(h_1-h_2)\charfun_\I+h_2\charfun$$
where $\charfun$ is the ``all $1$'s'' map. Recall from Lemma \ref{eigen} that
$h_1-h_2$ is an eigenvalue of $A$. So a \textit{positive} intriguing set has
$h_1-h_2=e^+$ and a \textit{negative} intriguing set has $h_1-h_2=e^-$. The
number of points of the partial quadrangle is
$$\frac{s(t+1)(\mu+st)}{\mu}+1.$$

%
%

\section{Intriguing sets of the known thin partial quadrangles}\label{thin}

A \textit{thin} partial quadrangle is simply a triangle-free strongly regular
graph. There are only seven known such graphs (see \cite[Chapter
8]{CameronvanLint}) and we explore and classify below the intriguing sets of
these geometries, for which many of the well-known interesting regular subgraphs
of these graphs predominate. Firstly, it is not difficult to see that the
pentagon contains no intriguing sets. The Petersen graph was dealt with in
Section \ref{alggraph}, and so it remains to consider the Clebsch,
Hoffman-Singleton, Gewirtz, $M_{22}$ and Higman-Sims graphs.

\subsection*{The Clebsch graph on 16 vertices}

In the Clebsch graph on $16$ vertices, the only negative intriguing sets are the
ten subgraphs isomorphic to $4K_2$, each stabilised by a group of order $192$;
which are maximal subgroups of the full group $2^4:S_5$.  As for positive
intriguing sets, the only examples are the forty $C_4$'s, a disjoint pair of
$C_4$'s, and complements of these. The Clebsch graph is small enough that we can
give a simple computer-free proof for the negative intriguing sets. Here is a
commonly used construction of the Clebsch graph. We have a special vertex
$\infty$, a set of five vertices $V_1=\{1,2,3,4,5\}$ and the subsets of $V_1$ of
size two, which we denote $V_2=\{12,13,14,15,23,24,25,34,35,45\}$. The vertex
$\infty$ is adjacent to all the members of $V_1$, the set $V_1$ is a coclique
and $V_2$ forms a Petersen graph whereby two elements are adjacent if they are
disjoint. Let $\mathcal{I}$ be a negative intriguing set of the Clebsch graph
with parameters $(h_1,h_2)$. Since the Clebsch graph is vertex transitive, we
may suppose that $\infty\in\mathcal{I}$. Moreover, the stabiliser of $\infty$
has $V_1$ and $V_2$ as two of its orbits, so we may also suppose without loss of
generality that $1\in\mathcal{I}$. Since $V_1$ is a coclique, there are no
further elements of $V_1$ inside $\mathcal{I}$, and we know now that $h_1=1$. In
fact, $\mathcal{I}$ must be a union of edges and have size $8$, as $h_2=4$. No
element of $V_2$ adjacent to $1$ can be in $\mathcal{I}$, so we can consider
$12$ and $15$ as external elements. For there to be $4$ elements in
$\mathcal{I}$ adjacent to $12$ (resp. $15$), we must have that $34$ and $45$ are
in $\mathcal{I}$.  The only edges of $V_2$ with no vertex adjacent to $1$ are
$\{34,25\}$, $\{45,23\}$ and $\{35,24\}$. By considering $15$, we see that all
of these edges must also be inside $\mathcal{I}$ and so it follows that
$\mathcal{I}=\{\infty,1,34,25,45,23,35,24\}$.  Hence the only negative
intriguing sets are the ten subgraphs isomorphic to $4K_2$.  Alternatively, we
can look to the minimal idempotent $E$ which annihilates $\charfun_\I$:
$$-\tfrac{1}{8}\left(\tiny{\begin{array}{rrrrrrrrrrrrrrrr}
 -5&   1&   1&   1&   1&  -1&   1&   1&   1&  -1&   1&   1&  -1&  -1&  -1&   1\\
  1&  -5&   1&   1&   1&   1&   1&  -1&   1&   1&   1&  -1&   1&  -1&  -1&  -1\\
  1&   1&  -5&   1&   1&   1&  -1&   1&   1&   1&  -1&   1&  -1&   1&  -1&  -1\\
   1&   1&   1&  -5&  -1&   1&   1&   1&  -1&   1&   1&   1&  -1&  -1&   1&  -1 \\
   1&   1&   1&  -1&  -5&   1&   1&   1&   1&  -1&  -1&  -1&   1&   1&  -1&   1 \\
  -1&   1&   1&   1&   1&  -5&   1&   1&  -1&   1&  -1&  -1&   1&   1&   1&  -1 \\
   1&   1&  -1&   1&   1&   1&  -5&   1&  -1&  -1&   1&  -1&   1&  -1&   1&   1 \\
   1&  -1&   1&   1&   1&   1&   1&  -5&  -1&  -1&  -1&   1&  -1&   1&   1&   1 \\
   1&   1&   1&  -1&   1&  -1&  -1&  -1&  -5&   1&   1&   1&   1&   1&  -1&   1 \\
  -1&   1&   1&   1&  -1&   1&  -1&  -1&   1&  -5&   1&   1&   1&   1&   1&  -1 \\
   1&   1&  -1&   1&  -1&  -1&   1&  -1&   1&   1&  -5&   1&   1&  -1&   1&   1 \\
   1&  -1&   1&   1&  -1&  -1&  -1&   1&   1&   1&   1&  -5&  -1&   1&   1&   1 \\
  -1&   1&  -1&  -1&   1&   1&   1&  -1&   1&   1&   1&  -1&  -5&   1&   1&   1 \\
   -1&  -1&   1&  -1&   1&   1&  -1&   1&   1&   1&  -1&   1&   1&  -5&   1&   1 \\
  -1&  -1&  -1&   1&  -1&   1&   1&   1&  -1&   1&   1&   1&   1&   1&  -5&   1 \\
    1&  -1&  -1&  -1&   1&  -1&   1&   1&   1&  -1&   1&   1&   1&   1&   1&  -5 
\end{array}}\right).$$
The points $\infty$ and $1$ in the above argument correspond to the first and
sixth rows above, which add to $-\tfrac{1}{8}( -6, 2, 2, 2, 2, -6, 2, 2, 0, 0,
0, 0, 0, 0, 0, 0 )$. The only way to cancel this vector out by adding other rows
of $E$, is to use all remaining rows or just the remainder from the first eight
rows.

\subsection*{The Hoffman-Singleton graph on 50 vertices}

The intriguing sets of the Hoffman-Singleton graph (on $50$ vertices) correspond
naturally to the maximal subgroups of its automorphism group
$\mathsf{PSU}(3,5).2$.  For the negative intriguing sets, we have one-hundred
$15$-cocliques (stabilised by $A_7$), the $252$ subgraphs isomorphic to $5C_5$
(each stabilised by a $5_+^{1+2}:8:2$) and pairs of disjoint $15$-cocliques
(each stabilised by an $M_{10}$).  The positive intriguing sets are also very
interesting: the 525 Petersen subgraphs (stabiliser: $2S_5.2$), a pair of
disjoint Petersen subgraphs (stabiliser: $D_{20}$) and three disjoint Petersen
subgraphs (stabiliser: $\mathsf{GL}(2,3):2$). The remaining intriguing sets are
complements of those above, and by computer, these are fully classified.

\subsection*{The Gewirtz graph on 56 vertices}

By computer, the only negative intriguing sets of the Gewirtz graph (on $56$
vertices) are the forty-two $16$-cocliques, the $105$ subgraphs isomorphic to
$6C_4$, the 420 Coxeter subgraphs (the graph on the antiflags of the Fano plane)
and complements of these.  The only positive intriguing sets are isomorphic to
the six regular subgraphs on $14$ vertices shown below, and those of greater
size obtained from a union of disjoint subgraphs or the complement of such a
subgraph. The six different types of positive intriguing sets of size $14$ form
single orbits under the automorphism group of the Gewirtz graph.

\begin{table}[H]
\centering
\begin{tabular}{c|c||c|c}
Subgraph&Aut. group&Subgraph&Aut. group\\
\hline\hline
\begin{minipage}[b]{0.15\linewidth}
\centering
\includegraphics[width=2.5cm]{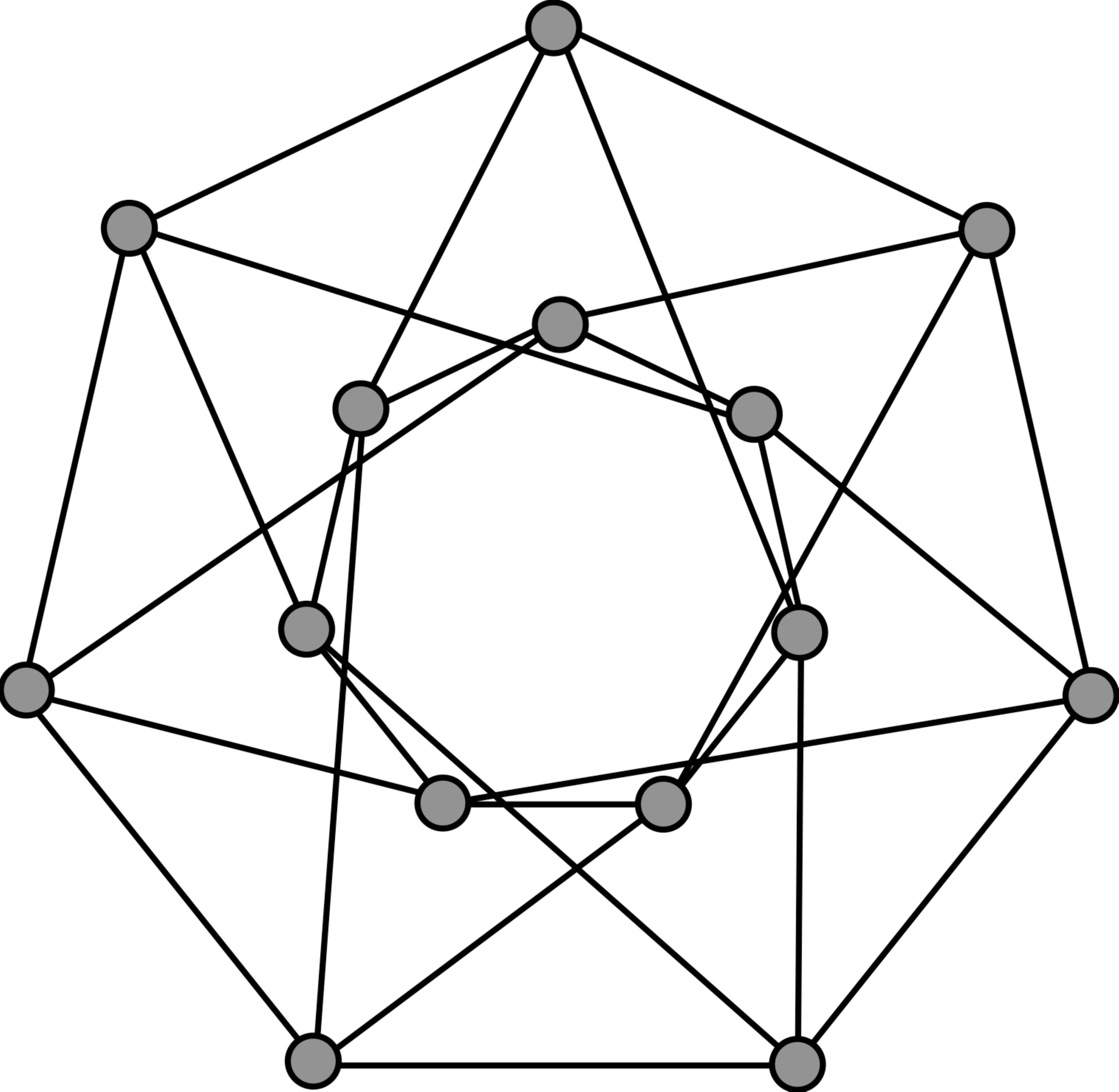}
\end{minipage}&\raisebox{6ex}{$D_{28}$}&
\begin{minipage}[b]{0.15\linewidth}
\centering
 \includegraphics[width=2.5cm]{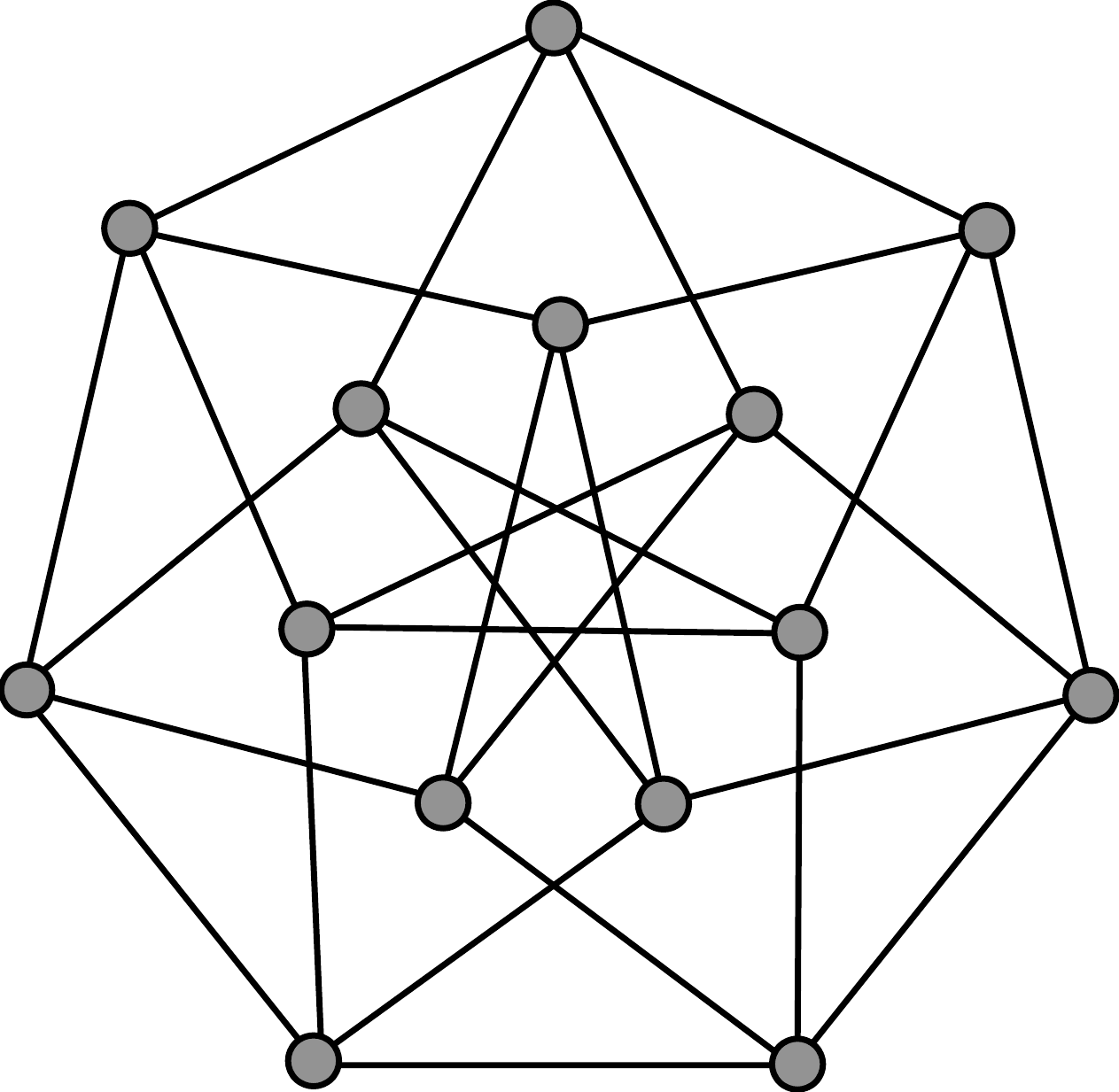}
\end{minipage}&\raisebox{6ex}{$D_{14}$}\\ \hline
\begin{minipage}[b]{0.15\linewidth}
\centering
 \includegraphics[width=2.5cm]{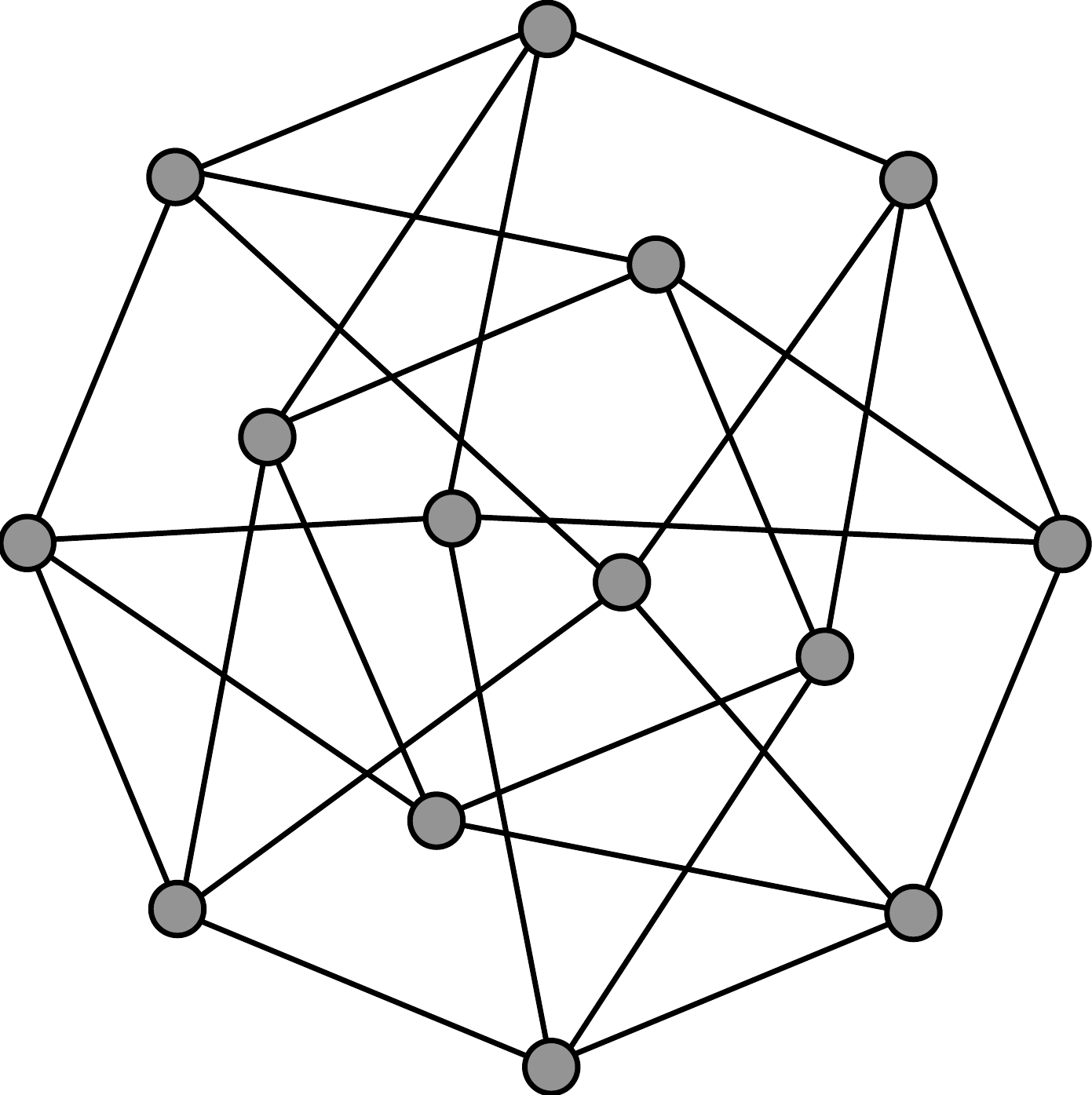}
\end{minipage}&\raisebox{6ex}{$D_8$}&
\begin{minipage}[b]{0.15\linewidth}
\centering
 \includegraphics[width=2.5cm]{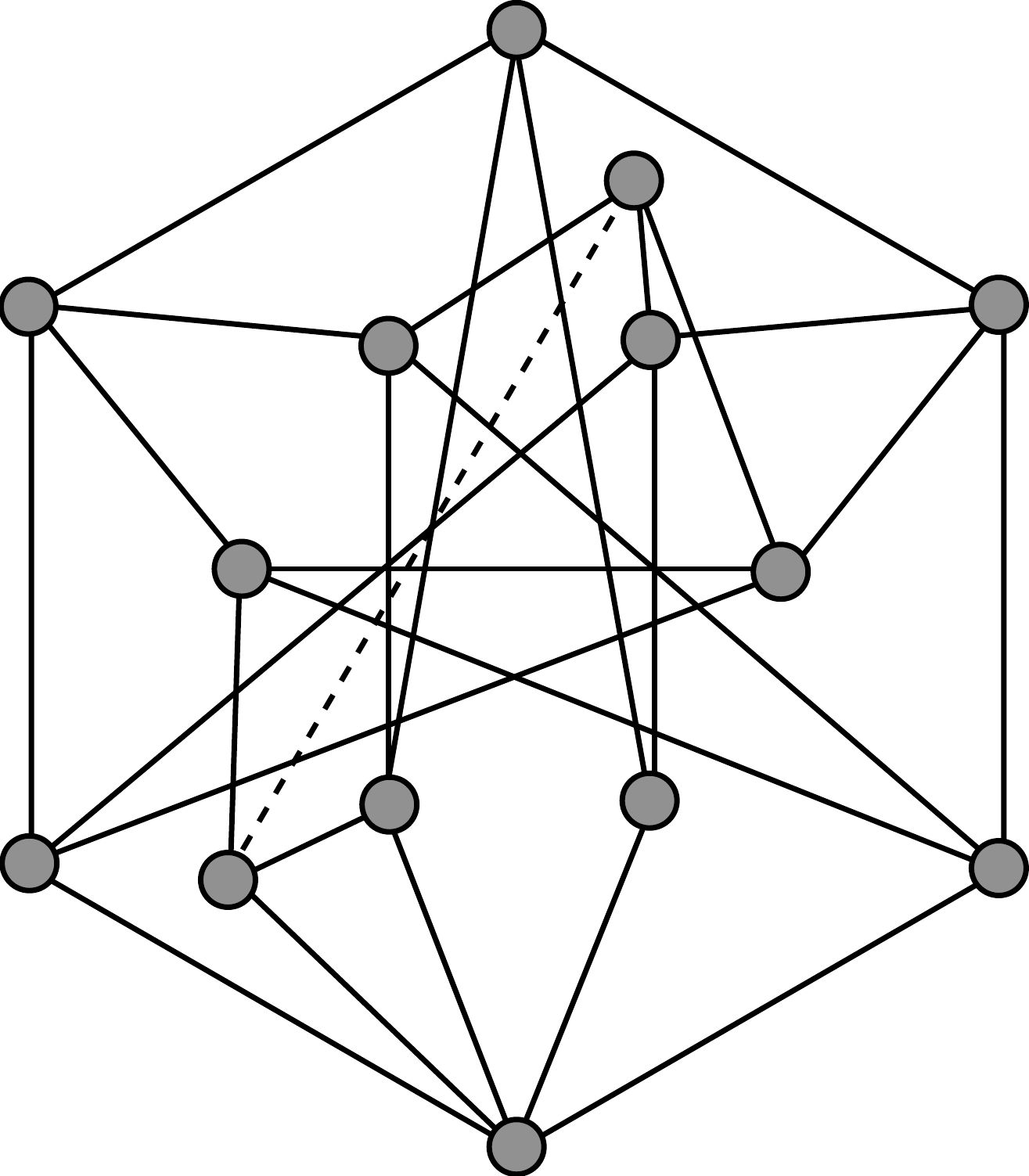}
\end{minipage}&\raisebox{6ex}{$D_{12}$}\\ \hline
\begin{minipage}[b]{0.15\linewidth}
\centering
 \includegraphics[width=2.5cm]{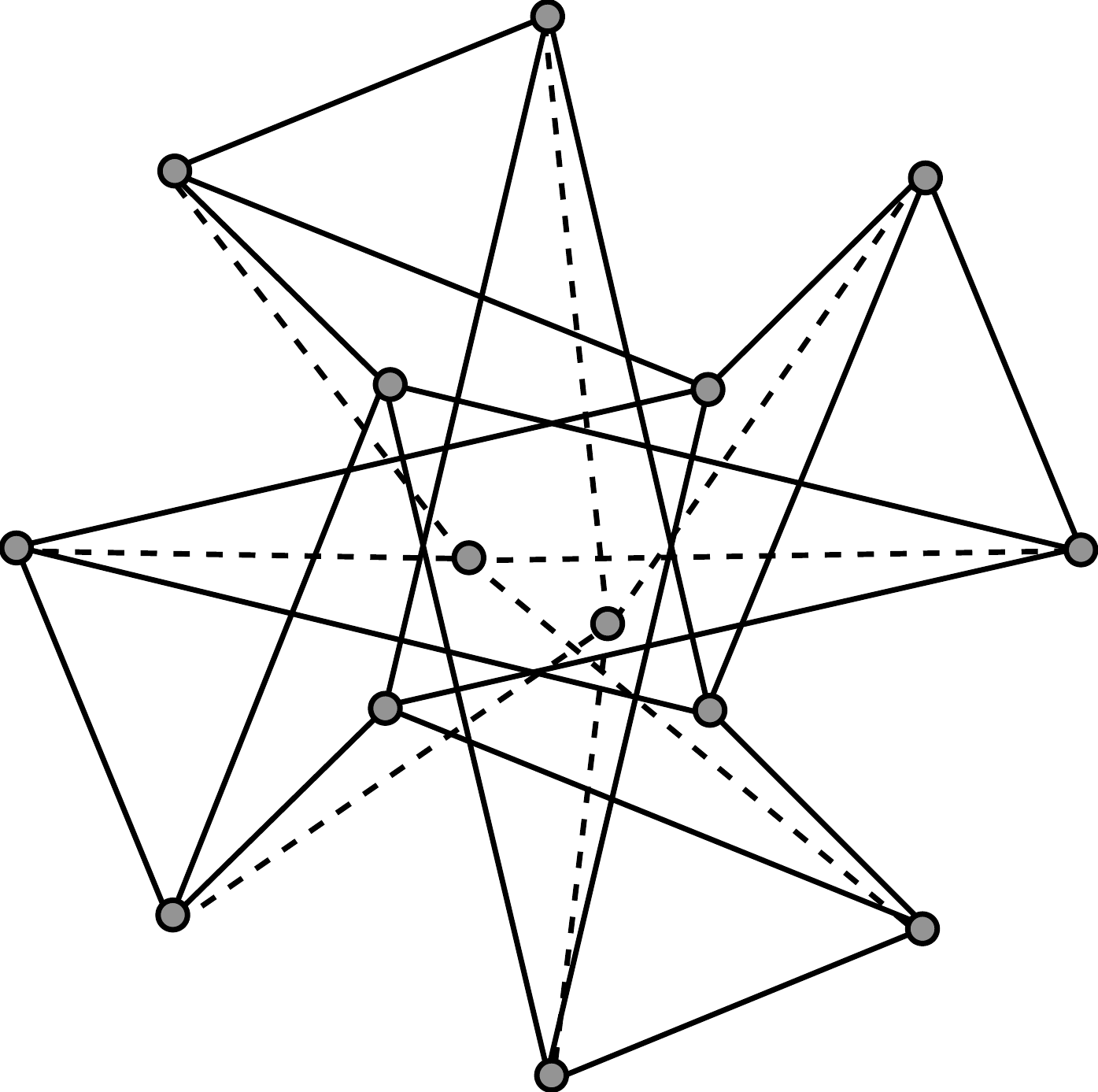}
\end{minipage}&\raisebox{6ex}{$C_2\times S_4$}&
\begin{minipage}[b]{0.15\linewidth}
\centering
 \includegraphics[width=2.5cm]{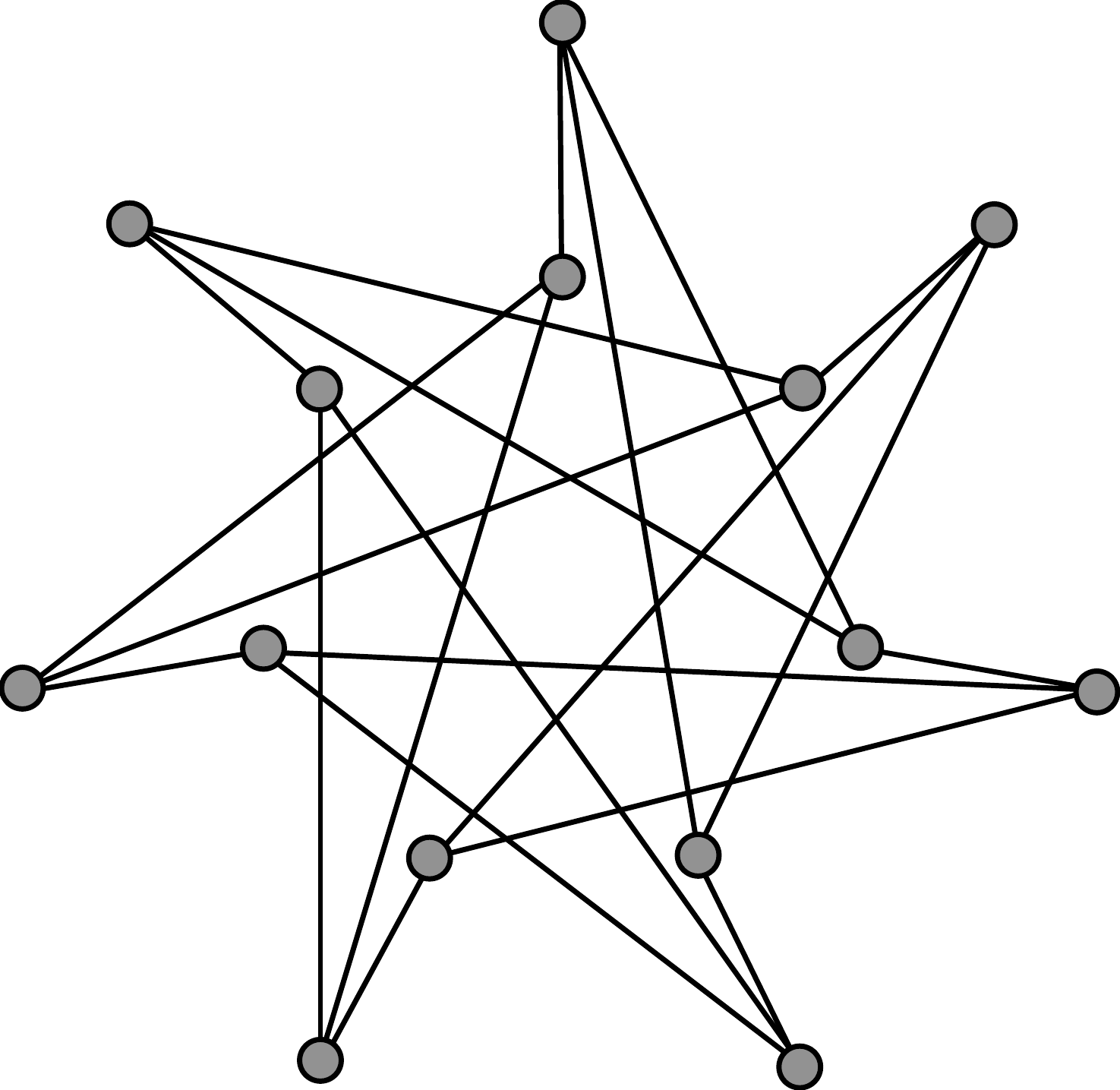}
\end{minipage}&\raisebox{6ex}{$\mathsf{PSL}(3,2):C_2$}
\end{tabular}
\caption{Positive intriguing sets of size $14$ in the Gewirtz graph. The first is a circulant and
  the last is the co-Heawood graph.}
\end{table}

\subsection*{The Higman-Sims $M_{22}$-graph on 77 vertices}

Two of the natural subgraphs of the $M_{22}$-graph are the $21$-cocliques and
the odd graphs $O_4$. These, and their complements, are the only negative
intriguing sets of the $M_{22}$-graph. A full classification of the positive
intriguing sets of the $M_{22}$-graph was not possible by computer, however, we
do have complete information of the positive intriguing sets which admit a
nontrivial automorphism group.  There are two interesting positive intriguing
sets which generate all the known examples. The first is a particular regular
subgraph on $11$ vertices (see the figure below) and the second is the incidence
graph of the complement of a biplane on $11$ points (i.e., $22$ vertices). There
exist disjoint triples of subgraphs of the first kind, and there exist disjoint
pairs consisting of one of each type of subgraph.

\begin{figure}[H]
\centering
 \includegraphics[scale=0.05]{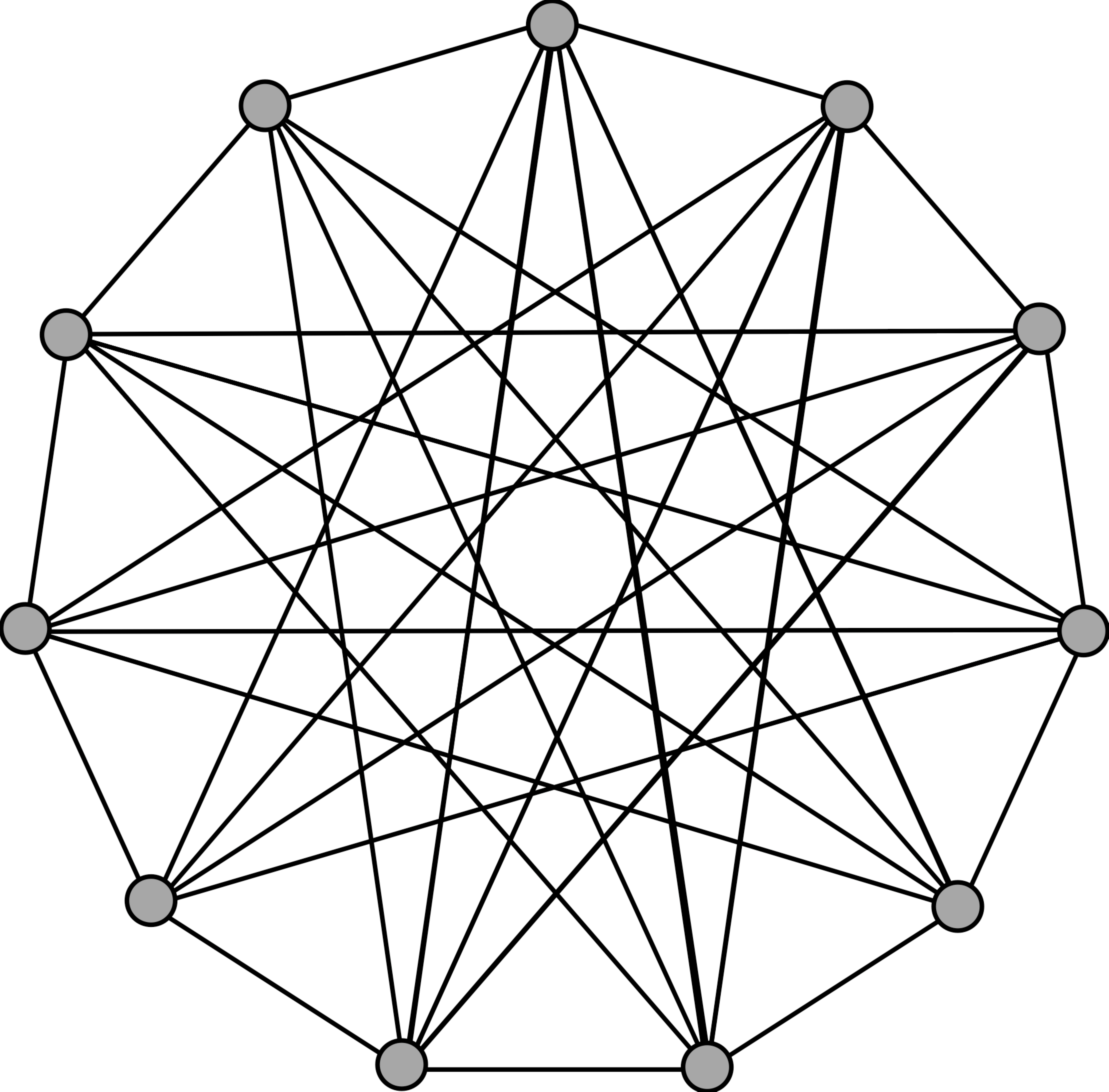}
 \caption{Circulant on $11$ vertices.}
\end{figure}

\subsection*{The Higman-Sims graph on 100 vertices}

The only negative intriguing sets of the Higman-Sims graph are the 704
Hoffman-Singleton subgraphs.  The known positive intriguing sets are as follows:
(i) a tetravalent circulant on $10$-vertices (see the figure below), (ii) the
graph which Brouwer \cite{Brouwer} calls $BD(K_5)$ (which is $K_{5,5}$ minus a
matching), (iii) bipartite on $20$ vertices, (iv) point-plane non-incidence
graph of $\mathsf{PG}(3,2)$ ($30$ vertices), (v) $2$-coclique extension of the
Petersen graph ($20$ vertices), (vi) a regular subgraph on $40$ vertices which
Brouwer \cite{Brouwer} denotes ``a pair of splits from the same family'', a
union of up to three disjoint subgraphs of type (i), and a union of up to three
disjoint subgraphs of type (ii).

 \begin{figure}[H]
 \centering
     \includegraphics[scale=0.04]{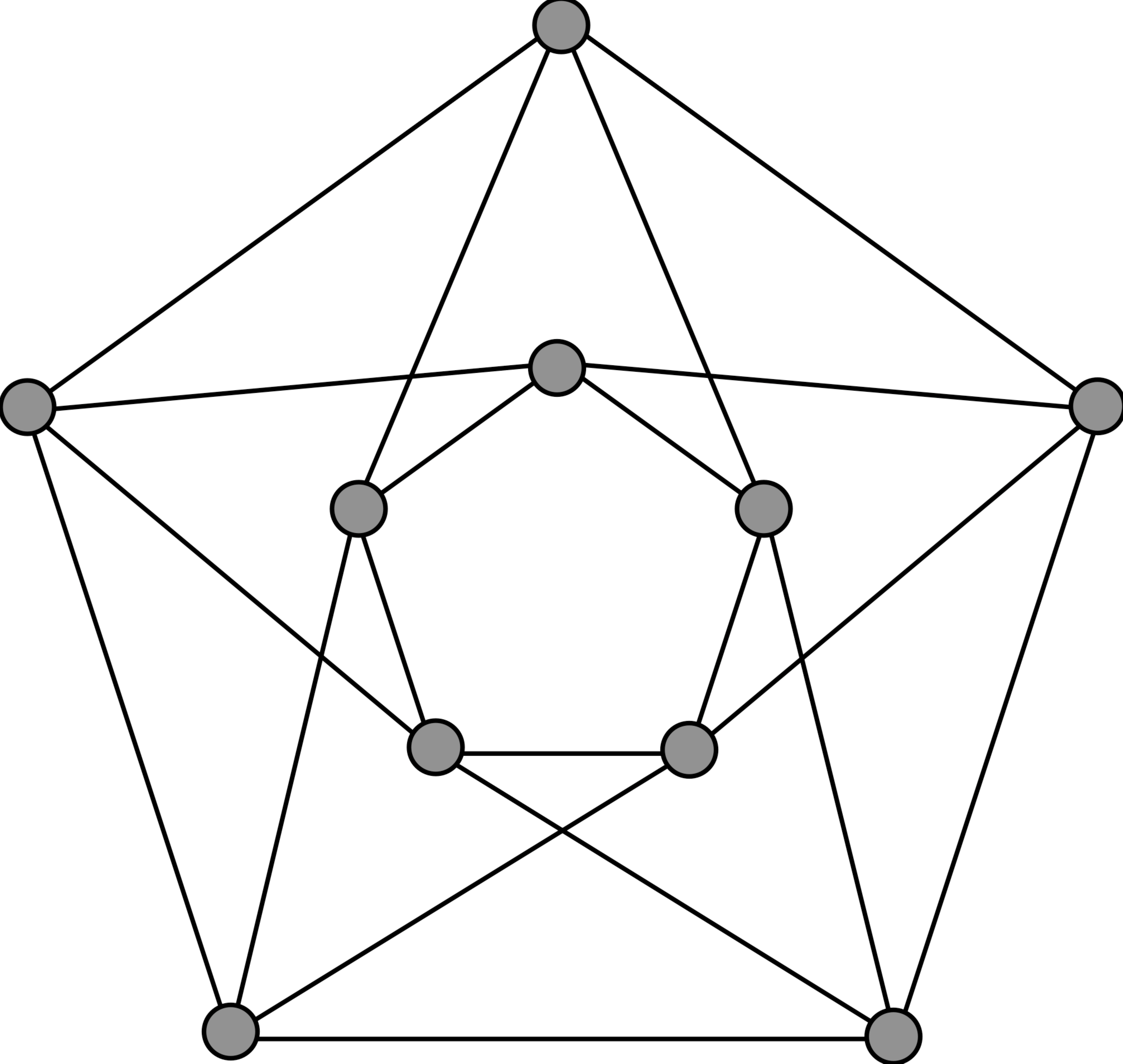}
 \caption{Tetravalent circulant on $10$ vertices.}
    \end{figure}


\section{Intriguing sets of generalised quadrangles and their interaction
with embedded partial quadrangles}\label{GQ}

Before we embark on an investigation into intriguing sets of partial quadrangles
which arise from point sets of generalised quadrangles, it will be necessary to
revise before-hand some of what we know about intriguing sets of generalised
quadrangles.  As was mentioned in the introduction, an intriguing set of a
generalised quadrangle is either an $m$-ovoid or an $i$-tight set. An $m$-ovoid
is a set of points such that every line meets it in $m$ points, and it is a
\textit{negative intriguing set} of the generalised quadrangle; that is, the
difference $h_1-h_2$ of its intersection numbers $h_1=m(t+1)-t-1$ and
$h_2=m(t+1)$ is negative (where $t+1$ is the number of lines on a point). An
$i$-tight set $\mathcal{T}$ is a set of points of a generalised quadrangle
$\mathcal{P}$ (of order $(s,t)$) such that the average number of points of
$\mathcal{T}$ collinear with a given point of $\mathcal{P}$ equals the maximum
possible value, namely $i+s$.  A set of points is \emph{tight} if it is
$i$-tight for some $i\ge 1$. The two intersection numbers here are $h_1= i+s-1$
and $h_2=i$, and so their difference $h_1-h_2$ is positive.

As the name suggests, $m$-ovoids are generalisations of \textit{ovoids}. An
ovoid of a generalised quadrangle is a set of points which partitions the lines,
that is, a $1$-ovoid.  The simplest tight sets are the $1$-tight sets, which one
can prove are the lines of a generalised quadrangle \cite{Payne87}. Hence the
point set covered by a partial spread (a set of disjoint lines) is a ubiquitous
example of a tight set of points. For more information on intriguing sets of
generalised quadrangles, we refer the reader to \cite{BLP}.

The partial quadrangles that we study in the following two sections are subsets
of points of generalised quadrangles, and hence, we will make use of the
following notion of ``intriguing at infinity''.

\begin{definition}[Intriguing at infinity]
  Let $\mathcal{G}$ be a generalised quadrangle of order $(s,s^2)$ and let
  $\infty$ be a set of points of $\mathcal{G}$ such that
  $\mathcal{G}\setminus\infty$ is a partial quadrangle. Then a set of points
  $\I$ of $\mathcal{G}$ is said to be \emph{intriguing at infinity} (with
  respect to $\infty$) if there are two constants $a_1$ and $a_2$ such that
$$
|y^\perp\cap\I\cap \infty|=
\begin{cases}
a_1&y\in \I\setminus \infty\\
a_2&y\notin \I \cup \infty.
\end{cases}$$
\end{definition}

\begin{theorem}\label{atinfinity}
  Let $\mathcal{G}$ be a generalised quadrangle of order $(s,s^2)$ and let
  $\infty$ be a set of points of $\mathcal{G}$ such that $\mathcal{G}\setminus
  \infty$ is a partial quadrangle . Let $\I$ be an intriguing set of
  $\mathcal{G}$ with parameters $(h_1,h_2)$. Then $\I\setminus \infty$ is an
  intriguing set of the partial quadrangle $\mathcal{G}\setminus \infty$ if and
  only if $\I$ is intriguing at infinity.
\end{theorem}

\begin{proof}
  Let $A$ be the adjacency matrix of the point graph of $\mathcal{G}$ and let
  $B$ be the adjacency matrix for the point graph of $\mathcal{G}\setminus
  \infty$.  Let $S$ be the matrix whose rows are indexed by the points of
  $\mathcal{G}$, and whose columns are indexed by $\mathcal{G}\setminus\infty$,
  such that the $(i,j)$-th entry of $S$ is equal to $1$ if the $i$-th point is
  equal to the $j$-th point of $\mathcal{G}\setminus\infty$, and $0$
  otherwise. Then
$$S^TAS = B\text{ and }S^TS=I.$$
When we write
$\charfun_\mathcal{H}^{\mathsf{PQ}}$ we mean the function $\charfun_\mathcal{H}$
restricted to the partial quadrangle.  By supposition, we have that
$$A\charfun_\I=(h_1-h_2)\charfun_{\I}+h_2\charfun.$$
Denote by $\infty'$ the complement of $\infty$. Note that $\I$ is intriguing at
infinity if and only if there exist non-negative integers $a_1$ and $a_2$ such
that
$$S^TA\charfun_{\I\cap\infty}=(a_1-a_2)\charfun_{\I\setminus\infty}^{\mathsf{PQ}}+
a_2\charfun_{\infty'}^{\mathsf{PQ}}.$$ On the other hand, $\I\setminus\infty$ is
intriguing in the partial quadrangle if and only if there exist non-negative
integers $h_1'$ and $h_2'$ such that
$$B\charfun_{\I\setminus\infty}^{\mathsf{PQ}}=(h_1'-h_2')\charfun_{\I\setminus\infty}^{\mathsf{PQ}}+
h_2'\charfun_{\infty'}^{\mathsf{PQ}}.$$

Now
$
A\charfun_{\I\cap\infty}=A(\charfun_\I-\charfun_{\I\cap\infty'})
=(h_1-h_2)\charfun_\I+h_2\charfun-A\charfun_{\I\setminus\infty}
$
and so $\I$ is intriguing at infinity if and only if there
exist non-negative integers $a_1$ and $a_2$ such that
$$
(a_1-a_2)\charfun^{\mathsf{PQ}}_{\I\setminus\infty}+a_2\charfun^{\mathsf{PQ}}_{\infty'}
=(h_1-h_2)\charfun^{\mathsf{PQ}}_{\I\setminus\infty}+h_2\charfun^{\mathsf{PQ}}_{\infty'}-S^TA\charfun_{\I\setminus\infty}.
$$
When we rearrange this equation, we obtain
$$
S^TA\charfun_{\I\setminus\infty}=
((h_1-a_1)-(h_2-a_2))\charfun^{\mathsf{PQ}}_{\I\setminus\infty}+(h_2-a_2)\charfun^{\mathsf{PQ}}_{\infty'}
$$
which is equivalent to $ B\charfun^{\mathsf{PQ}}_{\I\setminus \infty}=
((h_1-a_1)-(h_2-a_2))\charfun^{\mathsf{PQ}}_{\I\setminus\infty}+(h_2-a_2)\charfun^{\mathsf{PQ}}_{\infty'}$.
Therefore $\I\setminus \infty$ is an intriguing set of the partial quadrangle
$\mathcal{G}\setminus \infty$ if and only if $\I$ is intriguing at infinity.
\EndResult\end{proof}


\section{Partial quadrangles obtained by removing a point from a generalised quadrangle}\label{GQminusperp}

Let $\mathcal{G}$ be a generalised quadrangle of order $(s,t)$ and let $P$ be a
point of $\mathcal{G}$. Then the derived geometry with
\begin{center}
\begin{tabular}{l|l}
\textsc{points}& the points of $\mathcal{G}$ not collinear to $P$\\
\textsc{lines}& the lines of $\mathcal{G}$ not incident with $P$.
\end{tabular}
\end{center}

is a \textit{$(0,1)$-geometry}, that is, for every point $P$ and line $\ell$
which are not incident in this geometry, there is at most one point on $\ell$
collinear with $P$. The point graph of this geometry will be strongly regular if
and only if there is a constant $c$ such that for any two noncollinear points
$X$ and $Y$ of $\mathcal{G}$, not in $P^\perp$, there $c$ points of
$\mathcal{G}$ which are collinear with all three points $X$, $Y$ and $P$. This
property occurs when and only when the parameter $t$ is equal to $s^2$ (see
\cite{BoseShrikhande72} or \cite{Cameron75}), in which case $c=s+1$, and then we
obtain a partial quadrangle with parameters $(s-1,s^2,s(s-1))$.  In the
following lemma, we summarise the algebraic data needed to work with these kinds
of partial quadrangles.

\begin{lemma}\label{glossaryGQminusP}
  Let $\mathcal{G}$ be a generalised quadrangle of order $(s,s^2)$, let $P$ be a
  point of $\mathcal{G}$, and let $\I$ be an intriguing set of the partial
  quadrangle $\mathcal{G}\setminus P^\perp$ with intersection numbers
  $(h_1',h_2')$. Then we have the following information:

\begin{table}[H]
\begin{tabular}{l|c|c}
Case&Associated eigenvalue&Size\\
\hline
Negative intriguing set&$-s^2+s-1$&$h_2's$\\
Positive intriguing set&$s-1$&$h_2's^2/(s-1)$\\
Point set&$(s-1)(s^2+1)$&$s^4$\\
\hline
\end{tabular}
\caption{Eigenvalues and sizes of intriguing sets of $\mathcal{G}\setminus P^\perp$.}\label{table:glossaryGQminusP}
\end{table}
\end{lemma}

\begin{theorem}\label{atinfinity1}
  Let $\mathcal{G}$ be a generalised quadrangle of order $(s,s^2)$ and let
  $\infty=P^\perp$ where $P$ is a point of $\mathcal{G}$.  Let $\I$ be an
  intriguing set of $\mathcal{G}$ with parameters $(h_1,h_2)$ and which is
  intriguing at infinity with parameters $(a_1,a_2)$. Then $\I\setminus \infty$
  is an intriguing set with the same parity as $\I$ and we have the following
  possibilities for $(a_1,a_2)$:

\begin{table}[H]
\begin{tabular}{l|l|c|c|c}
Parity& Case &  $a_1$ & $a_2$ & $|\I\cap P^\perp|$\\
\hline
$m$-ovoid&$P\notin\I$& $m(s+1)-s$ & $m(s+1)$ &$m(s^2+1)$\\
&$P\in\I$& $m(s+1)-2s$ & $m(s+1)-s$ &$m(s^2+1)-s^2$\\
$i$-tight set&$P\notin\I$& $i/s$  & $i/s$ &$i$\\
&$P\in\I$& $(i-1)/s+1$ & $(i-1)/s+1$ &$i+s$
\end{tabular}
\caption{Possibilities for $(a_1,a_2)$.}
\end{table}
\end{theorem}

\begin{proof}
  Recall that the negative and positive eigenvalues for $\mathcal{G}$ are
  $-s^2-1$ and $s-1$, whilst they are $-s^2+s-1$ and $s-1$ for
  $\mathcal{G}\setminus P^\perp$. However, we must have that $h_1-h_2$ and
  $(h_1-a_1)-(h_2-a_2)$ are eigenvalues for the respective geometries.  In Table
  \ref{table:atinfinity1}, we outline the possibilities for these values
  depending on the four possible cases. We use the notation ``$-\rightarrow +$''
  (for example) to denote the case that $\I$ is negative intriguing and
  $\I\setminus\infty$ is positive intriguing.

\begin{table}[H]
\begin{center}
\begin{tabular}{l|c|c||c}
&$h_1-h_2$&$(h_1-a_1)-(h_2-a_2)$&$a_1-a_2$\\
\hline
$-\rightarrow -$&$-s^2-1$&$-s^2+s-1$&$-s$\\
$-\rightarrow +$&$-s^2-1$&$s-1$&$-s^2-s$\\
$+\rightarrow -$&$s-1$&$-s^2+s-1$&$s^2$\\
$+\rightarrow +$&$s-1$&$s-1$&$0$\\
\hline
\end{tabular}
\end{center}
\caption{Eigenvalues for the four possible cases.}\label{table:atinfinity1}
\end{table}

Since $a_1,a_2\le s^2+1$, we can rule out immediately the second case
above. Moreover, since $|Y^\perp \cap P^\perp \cap Z^\perp|=s+1$ for any three
pairwise non-collinear points $Y,P,Z$, and since there exists a point $Y\in
\I\setminus \infty$ and a point $Z \in (\mathcal{G}\setminus \infty)\setminus
\I$, it follows that $a_1-a_2\le s^2-s$. So the third case
in the above list is ruled out too. Hence parity is preserved. Now we see what
happens at infinity. Recall from Lemma \ref{eigen} that if $\I$ has associated
eigenvalue $e$ and $\I\setminus P^\perp$ has associated eigenvalue $e'$ (in the
partial quadrangle), then
$$|\I|=\frac{h_2}{s(s^2+1)-e}(s+1)(s^3+1)\quad\text{ and }\quad
|\I\setminus P^\perp|=\frac{h_2-a_2}{(s-1)(s^2+1)-e'}s^4.$$

Since $\I$ is intriguing we have
$$|P^\sim\cap\I|=\begin{cases}
h_2+e&P\in\I\\
h_2&P\notin\I.
\end{cases}$$

\noindent\textbf{Negative case:}
In the first case of Table \ref{table:atinfinity1}, $s^2+1$ divides $h_2$, and
$$|\I\cap P^\perp|=|\I|-|\I\setminus P^\perp|=\frac{h_2(s^3+1)}{s^2+1}-(h_2-a_2)s
=a_2s-\frac{h_2(s-1)}{s^2+1}.$$ As we know that $\I$ is an $m$-ovoid
(for some $m$), $h_2=m(s^2+1)$ and we have the following values:

\begin{center}
\begin{tabular}{l|c|c|c}
Case &  $a_1$ & $a_2$ & $|\I\cap P^\perp|$\\
\hline
$P\notin\I$& $m(s+1)-s$ & $m(s+1)$ &$m(s^2+1)$\\
$P\in\I$& $m(s+1)-2s$ & $m(s+1)-s$ &$m(s^2+1)-s^2$\\
\end{tabular}
\end{center}

\medskip

\noindent\textbf{Positive case:}
In the last case of Table \ref{table:atinfinity1}, we have
$$|\I\cap P^\perp|=|\I|-|\I\setminus P^\perp|=h_2(s+1)-\frac{h_2-a_2}{s-1}s^2=\frac{a_2s^2-h_2}{s-1}.$$
As we know that $\I$ is an $i$-tight set (for some $i$), $h_2=i$
and we have the following values:

\begin{center}
\begin{tabular}{l|c|c|c}
Case &  $a_1$ & $a_2$ & $|\I\cap P^\perp|$\\
\hline
$P\notin\I$& $i/s$  & $i/s$ &$i$\\
$P\in\I$& $(i-1)/s+1$ & $(i-1)/s+1$ &$i+s$\\
\end{tabular}
\end{center}

\EndResult\end{proof}


\subsection{Positive intriguing sets}

We now characterise the positive intriguing sets of a partial quadrangle
obtained from removing the perp of a point of a generalised quadrangle
$\mathcal{G}$, which are induced from intriguing sets of $\mathcal{G}$.

\begin{theorem}[Positive Intriguing $\longleftrightarrow$ Lines at Infinity]\label{positive}\ \\
  Let $\mathcal{G}$ be a generalised quadrangle of order $(s,s^2)$ and let $P$
  be a point of $\mathcal{G}$. If $\I$ is a positive intriguing set of
  $\mathcal{G}$, then $\I\setminus\infty$ is an intriguing set of
  $\mathcal{G}\setminus\infty$ if and only if $\I\cap \infty$ consists of $y$
  lines through $P$.
\end{theorem}

\begin{proof}
  First suppose that $\I$ is a positive intriguing set of $\mathcal{G}$. Then
  $|\I|=(s+1)i$, for some $i$. Assume that $\I$ intersects $\infty$ in $y$ lines
  through $P$. Then we have
$$
y=
\begin{cases}
i/s &P\notin \I\\
(i-1)/s+1 &P\in \I
\end{cases}$$ and it follows that $\I$ is intriguing at infinity with parameter
$y$, and hence by Lemma \ref{atinfinity}, $\I\setminus\infty$ is an intriguing
set of $\mathcal{G}\setminus\infty$.

Conversely, let $\I$ be a positive intriguing set of $\mathcal{G}$ and suppose
that $\I$ is intriguing at infinity with parameters $(a_1,a_2)$. By Lemma
\ref{atinfinity}, we have that $a_1=a_2=y$ with
$$
y=
\begin{cases}
i/s &P\notin \I\\
(i-1)/s+1 &P\in \I.
\end{cases}$$
By counting pairs $(Y,(Z,Z'))$ with $Y \in\I\setminus\infty$ and $Z,Z'\in \I\cap\infty$
with $Y \sim Z$, $Y\sim Z'$ and $Z\not\sim Z'$, we have
\[ s^4y(y-1)/2 = s^2 x\] where $x$ denotes the number of pairs $(Z,Z')$ as
above. Hence $x =s^2 y(y-1)/2 $. From Lemma \ref{atinfinity}, we also know that
the equation $|\I\cap\infty\setminus \{P\}|=ys$ is independent of whether $P\in
\I$ or $P\notin \I$. Finally it is easy to see that a set of $ys$ points in
$P^\perp\setminus \{P\}$ has the minimum number $y(y-1)/2 s^2$ of non-collinear
pairs (or the maximum number $ys(s-1)/2$ of collinear pairs) if and only if it
consists of $y$ lines through $P$. Indeed, we can see this by induction on
$y$. If $y=1$, then the statement is obvious. Assume it is true for a set of
$(y-1)s$ points. Take a set of $ys$ points. Consider any subset consisting of
$(y-1)s$ points. It has $(y-1)s(s-1)/2$ collinear pairs, otherwise we could not
have $ys(s-1)/2$ collinear pairs in total.  So by the induction hypothesis, we
have that the subset consists of $y-1$ lines through $P$. It now follows easily
that our point-set consists of $y$ lines through $P$, otherwise we would have
less collinear pairs.
\EndResult\end{proof}

The following also follow from Theorem \ref{positive}.

\begin{corollary}[Grid $\longrightarrow$ Positive Intriguing]\label{fromgrid}\ \\
  Let $\mathcal{G}$ be a generalised quadrangle of order $(s,s^2)$ and let $P$
  be a point of $\mathcal{G}$. For a point $X$ not collinear with $P$ and two
  distinct lines $\ell_1$, $\ell_2$ on $P$, let $[X,\ell_1,\ell_2]$ be the grid
  on $X$, $\ell_1$ and $\ell_2$. Then $[X,\ell_1,\ell_2]$ is a positive
  intriguing set of the partial quadrangle $\mathcal{G}\setminus P^\perp$, with
  parameters $(2s-2,s-1)$.
\end{corollary}

\begin{corollary}[$\mathsf{Q}(4,q)$ $\longrightarrow$ Positive Intriguing]\label{Q4q}\ \\
  Let $\mathcal{G}$ be the generalised quadrangle $\mathsf{Q}^-(5,q)$ and let
  $P$ be a point of $\mathcal{G}\setminus P^\perp$. Consider a $\mathsf{Q}(4,q)$
  embedded in $\mathsf{Q}^-(5,q)$. Then $\mathsf{Q}(4,q)\setminus P^\perp$ is a
  positive intriguing set of $\mathcal{G}\setminus P^\perp$ if and only if
  $\mathsf{Q}(4,q)\cap P^\perp$ is a tangent hyperplane to $\mathsf{Q}(4,q)$.
\end{corollary}


\subsection{Negative intriguing sets}

Segre \cite{Segre65} proved that if an $m$-ovoid of $\mathsf{Q}^-(5,q)$ exists,
then $m=(q+1)/2$; that is, it is a hemisystem. This fact can be readily extended
to all generalised quadrangles of order $(s,s^2)$, with an extra condition.

\begin{theorem}[Negative Intriguing $\longrightarrow$ Hemisystems]\label{negative}\ \\
  Let $\mathcal{G}$ be a generalised quadrangle of order $(s,s^2)$ and let $P$
  be a point of $\mathcal{G}$. If $\I$ is an $m$-ovoid of $\mathcal{G}$ and
  $\I\setminus\infty$ is an intriguing set of $\mathcal{G}\setminus\infty$, then
  $\I$ is a hemisystem $(m=(s+1)/2)$ of $\mathcal{G}$.
\end{theorem}

\begin{proof}
  Let $\I$ be an $m$-ovoid of $\mathcal{G}$, with $0<m<s+1$, and suppose that
  $\I$ is a negative intriguing set of $\mathcal{G}\setminus\infty$. Then by
  Lemma \ref{atinfinity} and Theorem \ref{atinfinity1}, $\I$ is intriguing at
  infinity with parameters $(a_1,a_2)$, where
$$
a_1=
\begin{cases}
m(s+1)-s &P\notin \I\\
m(s+1) -2s &P\in \I
\end{cases}$$
and
$$
a_2=
\begin{cases}
m(s+1) &P\notin \I\\
m(s+1)-s &P\in \I
\end{cases}.$$
Moreover
$$
|\I\cap\infty|=
\begin{cases}
m(s^2+1) &P\notin \I\\
m(s^2+1)-s^2 &P\in \I
\end{cases}.$$ First assume $P\in \I$, then $|\I\setminus\infty|=ms^3-ms^2+s^2$.
Counting pairs $(Y,(Z,Z'))$ with $Y \in \I\setminus\infty$ and $Z,Z'\in
\I\cap\infty$ with $Y \sim Z$, $Y\sim Z'$ and $Z\not\sim Z'$ we have
\[ a_1(a_1-1)(ms^3 -ms^2+s^2)/2 + a_2(a_2-1)(s^4-ms^3 +ms^2-s^2)/2 = s^2 x\]
where $x$ denotes the number of pairs $(Z,Z')$ as above. Hence $x= a_1(a_1-1)(ms
-m+1)/2 + a_2(a_2-1)(s^2-ms +m-1)/2$. On the other hand, since $\I$ is an
$m$-ovoid of $\mathcal{G}$ and $P\in \I$ we also know that $x = (m-1)^2
(s^2+1)s^2/2$.  According to Table \ref{table:atinfinity1} we have that
$a_1=m(s+1)-2s$ and $a_2=a_1+s$, in this case. Comparing the two values of $x$
obtained, we have
$$ 2m^2-3(s+1)m+(s+1)^2=0$$ from which it follows that 
$m=(s+1)/2$.

Next assume that $P\notin \I$, then $|\I\setminus\infty|=ms^3-ms^2$.  Counting
pairs $(Y,(Z,Z'))$ with $Y \in \I\setminus\infty$ and $Z,Z'\in \I\cap\infty$
with $Y \sim Z$, $Y\sim Z'$ and $Z\not\sim Z'$ we have
\[ a_1(a_1-1) (ms^3 -ms^2)/2 + a_2(a_2-1)(s^4-ms^3 +ms^2)/2 = s^2 x\] where $x$
denotes the number of pairs $(Z,Z')$ as above. Hence $x= a_1(a_1-1)(ms -m)/2 +
a_2(a_2-1)(s^2-ms +m)/2$. On the other hand, since $\I$ is an $m$-ovoid of
$\mathcal{G}$ and $P\notin \I$ we also know that $x = m^2 (s^2+1)s^2/2$.
Comparing the two values of $x$ obtained, we have $$m(2m - (s+1))=0$$ from which
it follows that $m=(s+1)/2$.
\EndResult\end{proof}

Moreover, for hemisystems we have

\begin{lemma}[Hemisystem $\longrightarrow$
  Negative Intriguing]\label{hemi_negint}\label{NiceHemi}
  Let $\I$ be a hemisystem of a generalised quadrangle $\mathcal{G}$ of order
  $(s,s^2)$. Let $P$ be a point of $\mathcal{G}$. Then $\I\setminus P^\perp$ is
  a negative intriguing set of $\mathcal{G}\setminus P^\perp$ if and only if
  $|Y^\perp \cap \I\cap P^\perp|$ is a constant over all $P$ and $Y$ not both in
  $\I$.
\end{lemma}

\begin{proof}
Follows from Theorem \ref{atinfinity1}.
\EndResult\end{proof}

\medskip\noindent\textit{Open question:}
  For every $m$-ovoid $\mathcal{O}$ of a generalised quadrangle $\mathcal{G}$ of
  order $(s,s^2)$, does there exist a point $P$ such that $\mathcal{O}\setminus
  P^\perp$ is an intriguing set of the associated partial quadrangle?  (Compare
  with Theorem \ref{negative}).
\medskip

\begin{lemma}[Cone$\longrightarrow$ Negative Intriguing]\label{NiceCone}
  Let $\mathcal{G}$ be the generalised quadrangle of order $(s,s^2)$ and let $P$
  be a point of $\mathcal{G}$. For every point $Z\in P^\perp$, the set of points
  $Z^\perp\setminus P^\perp$ is a negative intriguing set of
  $\mathcal{G}\setminus P^\perp$ with parameters $(s-1,s^2)$.
\end{lemma}

\begin{proof}
  Let $Z$ be a point of $P^\perp$ and let $\I$ be the set of points of
  $\mathcal{G}\setminus P^\perp$ contained in $Z^\perp$.  (Clearly if $Z=P$ we
  get the empty set, so assume that $Z\ne P$). Let $X$ be a point in $\I$. Then
  the only points of $Z^\perp$ collinear with $X$ lie on the line
  $XZ$. Moreover, every point but $Z$ on this line is not in $P^\perp$. So there
  are $s-1$ points collinear with $X$ (and not equal to $X$) in $\I$. Now let
  $Y$ be a point not in $\I$, but in $\mathcal{G}\setminus P^\perp$. Now $Y$ is
  not collinear to $Z$, and we want to know how many points of $\mathcal{G}$ are
  collinear with both $Y$ and $Z$, but not on the line $ZP$.  This number is
  $\mu - 1 = s^2$. Therefore, $\I$ is a negative intriguing set of
  $\mathcal{G}\setminus P^\perp$ with parameters $(s-1,s^2)$.
\EndResult\end{proof}

So from Lemma \ref{NiceHemi} and Lemma \ref{NiceCone}, we have two ways to
obtain negative intriguing sets of a partial quadrangle which is a generalised
quadrangle minus the perp of a point: namely, from unions of disjoint cones, and
from hemisystems. We conjecture that these are the only two possible types of
negative intriguing sets.

\begin{conjecture}\label{negint_minusperp}\ \\
  Let $\mathcal{G}$ be a generalised quadrangle of order $(s,s^2)$ and let $P$
  be a point of $\mathcal{G}$. If $\I$ is a negative intriguing set of the
  partial quadrangle $\mathcal{G}\setminus P^\perp$, then either:
\begin{enumerate}
\item[(i)] There exists points $Z_1,\ldots,Z_n$ of $P^\perp$ such that
$\I=\bigcup_{i=1}^n \left(Z_i^\perp\setminus P^\perp\right),$ or
\item[(ii)] There exists a hemisystem $\mathcal{H}$ of $\mathcal{G}$ such that
$\I=\mathcal{H}\setminus P^\perp.$
\end{enumerate}
\end{conjecture}

We are able to provide a partial answer to the above conjecture via the lemma
and theorem below.  For the identity and ``all-ones'' matrices, we will
sometimes use a subscript which describes the size of the matrix. For example,
$I_{P^\perp}$ and $J_{P^\perp}$ denote the corresponding square matrices with
$|P^\perp|$ rows and columns.

\begin{lemma}\label{icky}
  Let $\mathcal{G}$ be a generalised quadrangle of order $(s,s^2)$ and let $P$
  be a point of $\mathcal{G}$.  Order the points of $\mathcal{G}$ so that the
  points of $P^\perp$ appear last, with $P$ last of all. Let $A$ be the
  adjacency matrix of the point-graph of $\mathcal{G}$, and let $B$ be the
  adjacency matrix of the partial quadrangle $\mathcal{G}\setminus P^\perp$ such
  that $$A=\begin{pmatrix}B&C\\ C^T& D\end{pmatrix}.$$ Let
  $\lambda=-s^2-1$. Then:
\begin{enumerate}
\item[(a)] $D-\lambda I_{P^\perp}$ is invertible and moreover
$$s^3(s^2+1)(s+1)(D-\lambda I)^{-1}=(s^4+s^3+s-1)I+J-(s^2+1) D-s(M+M^T)+s(s^2+s-1)E$$
where
$$M:=\left(
\begin{smallmatrix}
0&\cdots&0&0\\
\vdots&&\vdots&\vdots\\
0&\cdots&0&0\\
1&\ldots&1&0
\end{smallmatrix}
\right)\text{ and }
E:=\left(
\begin{smallmatrix}
0&\cdots&&0\\
\vdots&&&\vdots\\
&&&0\\
0&\cdots&0&1
\end{smallmatrix}
\right).$$

\item[(b)] $C(D-\lambda I)^{-1}\charfun_{P^\perp}=\charfun_{PQ}$.
\item[(c)] If $\I$ is a negative intriguing set of $\mathcal{G}\setminus
  P^\perp$ with parameters $(h_1',h_2')$, then
$$CC^T\charfun_\I=s^3\charfun_\I+s|\I|\charfun_{PQ}$$ and
$$C(D-\lambda I)^{-1}C^T\charfun_{\I}=s\charfun_\I+h_2'\charfun_{PQ}.$$
\end{enumerate}
\end{lemma}

\begin{proof}
  \medskip\noindent\textbf{(a)} Since $D$ is the adjacency matrix for $P^\perp$,
  the eigenvalues for $D$ are $s-1$ and $-1$. Since $\lambda < -1$, it follows
  that $D-\lambda I_{P^\perp}$ is invertible.  Now we apply $D-\lambda
  I_{P^\perp}$ to our proposed formula for the inverse $(D-\lambda
  I_{P^\perp})^{-1}$:
\begin{align*}
&(D-\lambda I_{P^\perp})((s^4+s^3+s-1)I+J-(s^2+1) D-s(M+M^T)+s(s^2+s-1)E)=\\
&(s^4+s^3+s-1)D+DJ-(s^2+1) D^2-s(DM+DM^T)+s(s^2+s-1)DE-\\
&\lambda(s^4+s^3+s-1)I-\lambda J+\lambda(s^2+1) D+\lambda s(M+M^T)-\lambda s(s^2+s-1)E.
\end{align*}

Recall that the last row and column of $D$ represent the point $P$.  To compute
the $(i,j)$-entry of $D^2$, we note that if $i,j\ne s^3+s+1$, then
$$D^2(i,j)=J(i,j)+(s-1)I(i,j)+(s-2)D(i,j).$$
So to complete the equation, we consider what happens when $i=j=s^3+s+1$. We
then see that
$$D^2=J+(s-1)I+(s-2)D+s^3E.$$
Also, it is not difficult to see that $DM=J-M-M^T-E$, $DE=M^T$, $DJ=sJ+s^3(M+E)$
and $DM^T=(s-1)M^T+(s^3+s)E$.  So our equation simplifies to the
following:
$$(D-\lambda I_{P^\perp})((s^4+s^3+s-1)I+J-(s^2+1) D-s(M+M^T)+s(s^2+s-1)E)=s^3(s^2+1)(s+1)I$$
from which the desired conclusion follows.

\medskip\noindent\textbf{(b)} Note that
$C\charfun_{P^\perp}=(s^2+1)\charfun_{PQ}$,
$CJ\charfun_{P^\perp}=(s^2+1)(s^3+s+1)\charfun_{PQ}$,
$CM\charfun_{P^\perp}=CE\charfun_{P^\perp}=\underline{0}$ (the zero vector) and
$CM^T\charfun_{P^\perp}=(s^2+1)\charfun_{PQ}$.
Now $$D\charfun_{P^\perp}=s\charfun_{P^\perp}+(0,\ldots,0,s^3)$$ and
hence $$CD\charfun_{P^\perp}=s(s^2+1)\charfun_{P^\perp}+C(0,\ldots,0,s^3)=s(s^2+1)\charfun_{PQ}.$$
A little calculation then shows that $C(D-\lambda
I)^{-1}\charfun_{P^\perp}=\charfun_{PQ}$.

\medskip\noindent\textbf{(c)} First we prove that
$CC^T\charfun_\I=s^3\charfun_\I+s|\I|\charfun_{PQ}$.  Let $P_i$ and $P_j$ be the
$i$-th and $j$-th points of the partial quadrangle.  The $(i,j)$ entry of $CC^T$
is the number of points of $P^\perp$ which are collinear with both $P_i$ and
$P_j$. Now if $P_i$ and $P_j$ are collinear, that is $B(i,j)=1$, then the only
point of $P^\perp$ collinear to both $P_i$ and $P_j$ is the point of
intersection of $P^\perp$ with the line joining $P_i$ and $P_j$; so
$CC^T(i,j)=1$ in this case. Otherwise, if $B(i,j)=0$, then there are $s+1$
points of $P^\perp$ collinear to all three points $P$, $P_i$ and $P_j$ (recall
that this was a property of the ambient generalised quadrangle $\mathcal{G}$ for
$\mathcal{G}\setminus P^\perp$ to be a partial quadrangle). Therefore,
$CC^T=(s+1)J_{PQ}-sB+(s^2-s)I_{PQ}$ and hence
\begin{align*}
  CC^T\charfun_\I&=(s+1)J_{PQ}\charfun_\I-sB\charfun_\I+(s^2-s)I_{PQ}\charfun_\I\\
  &=(s+1)|\I|\charfun_{PQ}-s((-s^2+s-1)\charfun_{\I}  +h_2'\charfun_{PQ})+(s^2-s)\charfun_\I\\
  &=s^3\charfun_\I+s|\I|\charfun_{PQ}.
\end{align*}

Now we list some formulae which can be worked out with some simple geometric
arguments:
\begin{align*}
CMC^T&=CM^TC^T=CEC^T=0,\\
CJC^T&=(s^2+1)^2J_{PQ},\\
CDC^T=(s^2+1)J_{PQ}-CC^T.
\end{align*}
The last of these formulae will serve as a demonstration of how to compute all
of them.  The matrix $DC^T$ measures the number of points of $P^\perp$ which are
collinear with two points, one from $P^\perp$ and the other from the partial
quadrangle.  Upon applying $C$, we see see that $CDC^T=(s^2+1)J_{PQ}-CC^T$.
From the above calculations, we arrive at
\begin{align*}
s^3(s^2+1)(s+1)C(D-\lambda I)^{-1}C^T\charfun_{\I}&=s(s^2+1)(s+1)CC^T\charfun_{\I}\\
&=s^3(s^2+1)(s+1)\left(s\charfun_\I+h_2'\charfun_{PQ}\right).
\end{align*}
Therefore, $C(D-\lambda I)^{-1}C^T\charfun_{\I}=s\charfun_\I+h_2'\charfun_{PQ}$.
\EndResult\end{proof}

\begin{theorem}[Negative intriguing set of the right size $\longleftrightarrow$ Hemisystem]
  Let $\mathcal{G}$ be a generalised quadrangle of order $(s,s^2)$, $s$ odd, and
  let $P$ be a point of $\mathcal{G}$ such that $\mathcal{G}\setminus P^\perp$
  is a partial quadrangle.  Suppose $\I$ is a negative intriguing set of the
  partial quadrangle $\mathcal{G}\setminus P^\perp$ such that $|\I|$ is either
$$s^2(s^2-1)/2\quad \text{  or  } \quad s^2(s^2+1)/2.$$
Then there is a subset $\I^*$ of points of $P^\perp$ such that $\I\cup
\I^*$ is a hemisystem of $\mathcal{G}$.
\end{theorem}

\begin{proof}
  Suppose that $\I$ is a negative intriguing set of $\mathcal{G}\setminus
  P^\perp$ with parameters $(h_1',h_2')$.  Let $\lambda$ be the negative
  eigenvalue of $A$ (i.e., $-s^2-1$), let $h_2$ be a positive integer, and
  let $$v=(D-\lambda I_{P^\perp})^{-1}(-C^T\charfun_I+h_2\charfun_{P^\perp}).$$
  We will show that there is a value of $h_2$ such that $v$ is a $(0,1)$-vector
  and hence represents a subset $\I^*$ of points of $P^\perp$. If we also show
  that $v$ corresponds naturally to an eigenvector of $A$ (see Lemma
  \ref{eigen}), then it will follow that $\I\cup \I^*$ is intriguing in
  $\mathcal{G}$. We show first that $\charfun_\I+v-\alpha \charfun_{GQ}$ is an
  eigenvector of $A$ with eigenvalue $\lambda$,
  where $$\alpha=\frac{h_2}{(s+1)(s^2+1)}.$$ We apply $A$ to our proposed
  eigenvector:
\begin{align*}
A\left(\charfun_\I+v-\alpha \charfun_{GQ}\right)=&A\charfun_\I+Av-\alpha s(s^2+1)\charfun_{GQ}\\
=&\begin{pmatrix}
B\charfun_\I+Cv\\ 
C^T\charfun_\I+Dv\end{pmatrix}-\alpha s(s^2+1)\charfun_{GQ}.
\end{align*}

By Lemma \ref{icky}, 
\begin{align*}
Cv&=-C(D-\lambda I)^{-1}C^T \charfun_{\I}+h_2\charfun_{PQ}=-s\charfun_\I-h_2'\charfun_{PQ}+h_2\charfun_{PQ}\\
&=-\left((\lambda+s)\charfun_\I+h_2'\charfun_{PQ}\right)+\lambda\charfun_\I+h_2\charfun_{PQ}
=-B\charfun_\I+\lambda\charfun_\I+h_2\charfun_{PQ}
\end{align*}

and hence $B\charfun_\I+Cv=\lambda\charfun_\I+h_2\charfun_{PQ}$. We also have
$$C^T\charfun_\I+Dv=
C^T\charfun_\I+\lambda v -C^T\charfun_{\I}+h_2\charfun_{P^\perp} =\lambda v +h_2\charfun_{P^\perp}.$$

Hence,
\begin{align*}
A(\charfun_\I+v-\alpha \charfun_{GQ})=&A\charfun_\I+Av-\alpha s(s^2+1)\charfun_{GQ}
=\begin{pmatrix}
B\charfun_\I+Cv\\  
C^T\charfun_\I+Dv\end{pmatrix}-\alpha s(s^2+1)\charfun_{GQ}\\
=&\begin{pmatrix}
\lambda(\charfun_\I-\alpha(s+1) \charfun_{PQ})\\  
\lambda (v -\alpha(s+1)\charfun_{P^\perp})\end{pmatrix}+\alpha s\lambda\charfun_{GQ}\\
=&\lambda (\charfun_\I+v-\alpha \charfun_{GQ}).
\end{align*}

So $\charfun_\I+v-\alpha \charfun_{GQ}$ is an eigenvector of $A$ with eigenvalue
$\lambda$.  Note that this is true no matter what choice we make for the value
of $h_2$. We will show that there exists a specific $h_2$ such that $v$ is a
$(0,1)$-vector. A tedious calculation shows that if we restrict $v$ to the
points of $P^\perp\setminus\{P\}$, then
$$v|_{P^\perp\setminus\{P\}}=\frac{-C^T\charfun_\I}{s^2}+\frac{h_2}{s^2+1}\charfun_{P^\perp\setminus\{P\}}.$$
We then look at the action of $v$ on the point $P$ to derive
$$v|_{\{P\}}=\frac{|\I|}{s^2}+\frac{h_2}{s^2+1}(1-s).$$
Let $X=h_2/(s^2+1)$. Then $\sum v=X(s^3+1)-|\I|$ and

\begin{align*}
v\cdot v&=\frac{\charfun_\I CC^T\charfun_\I^T}{s^4}-2X\frac{|\I|(s^2+1)}{s^2}+X^2(s^3+s)
+\left(\frac{|\I|}{s^2}+X(1-s)\right)^2\\
&=\frac{|\I|}{s}+\frac{(s+1)|\I|^2}{s^4}-2X\frac{|\I|(s+1)}{s}+X^2(s^3+s^2-s+1).
\end{align*}

So
\begin{align*}
  v\cdot v-\sum v =&\frac{|\I|}{s}+\frac{(s+1)|\I|^2}{s^4}-2X\frac{|\I|(s+1)}{s}+X^2(s^3+s^2-s+1)-X(s^3+1)+|\I|\\
  =&X^2(s^3+s^2-s+1)+X\left(\frac{-2|\I|(s+1)}{s}-s^3-1\right)+
  \frac{|\I|}{s}+\frac{(s+1)|\I|^2}{s^4}+|\I|.
\end{align*}

Now suppose $|\I|=s^2(s^2-1)/2$. Then it turns out that
$$v\cdot v-\sum v=\left((s^3+s^2-s+1)X-\frac{s+1}{2}(s^3+s^2-s+1)\right)\left(X-\frac{(s-1)(s^3+3s^2+s-1)}{2(s^3+s^2-s+1)}\right)$$
and hence $v$ is a zero vector if and only if $X=(s+1)/2$.  That is, if we let
$h_2=(s^2+1)(s+1)/2$, then $v$ is a characteristic function for a subset $\I^*$
of $P^\perp$, and the union of $\I$ with $\I^*$ forms a negative intriguing set
of the generalised quadrangle $\mathcal{G}$. By \cite[Theorem 4.1]{BLP},
$\I\cup\I^*$ is a hemisystem of $\mathcal{G}$ (as $h_2=(s^2+1)(s+1)/2$).  A
similar argument holds for the case $|\I|=s^2(s^2+1)/2$.
\EndResult\end{proof}


\subsection{Examples of intriguing sets of $\mathsf{Q}^-(5,q)$}

Segre \cite{Segre65} proved that for $q=3$, there is just one hemisystem up to
equivalence, and it was long thought to be the only example of such an
object. However, Cossidente and Penttila \cite{CossidentePenttila05} constructed
an infinite family of hemisystems of $\mathsf{Q}^-(5,q)$ admitting
$\mathsf{P\Omega}^-(4,q)$, together with a special example for $q=5$ admitting
the triple cover of $A_7$.

There are many tight sets of $\mathsf{Q}^-(5,q)$, simply because there are many
partial spreads of $\mathsf{Q}^-(5,q)$. However, some interesting examples arise
from Cameron-Liebler line classes. A set of lines $\mathcal{L}$ of
$\mathsf{PG}(3,q)$ is said to be a \textit{Cameron-Liebler line class} if there
exists a constant $i$ such that $\mathcal{L}$ meets every (regular) line spread
of $\mathsf{PG}(3,q)$ in $i$ elements. Such a set of lines gives rise to an
$i$-tight set of $\mathsf{Q}^-(5,q)$ as follows: first note that every spread of
the symplectic generalised quadrangle $\mathsf{W}(3,q)$ is a spread of
$\mathsf{PG}(3,q)$, and so the set of lines of $\mathcal{L}$ in
$\mathsf{W}(3,q)$ meets each spread of $\mathsf{W}(3,q)$ in $i$ elements. Hence,
by dualising, we obtain an $i$-tight set of $\mathsf{Q}(4,q)$. By embedding, we
produce an $i$-tight set of $\mathsf{Q}^-(5,q)$.  For
$i\in\{0,1,2,q^2-1,q^2,q^2+1\}$, it was shown by Cameron and Liebler
\cite{CameronLiebler82} that a Cameron-Liebler line class of $\mathsf{PG}(3,q)$
is one of the following:
\begin{enumerate}
\item[(i)] the empty set $(i=0)$ or its complement $(i=q^2+1)$;
\item[(ii)] the set of lines on a point ($i=1$) or its complement $(i=q^2)$;
\item[(iii)] the set of lines in a hyperplane ($i=1$) or its complement
  $(i=q^2)$;
\item[(iv)] the set of lines on a point $P$ together with the lines in a
  hyperplane $H$, where $P$ is not in $H$ ($i=2$), or its complement
  $(i=q^2-1)$.
\end{enumerate}
The Cameron-Liebler line classes above can only give rise to tight sets of
$\mathsf{Q}^-(5,q)$ which consist of a line, a pair of skew lines, or a
complement of one of these.  However, much more is known about the existence and
non-existence of Cameron-Liebler line classes, and so we refer the interested
reader to \cite{GovaertsPenttila} for more on this topic. Finally we note that a
$\mathsf{Q}(4,q)$ embedded in $\mathsf{Q}^-(5,q^2)$ (subfield embedding) is
$(q+1)$-tight, and the points of $\mathsf{Q}(4,q^2)$ which are collinear but not
equal to their conjugate forms a $q(q^2-1)$-tight set of $\mathsf{Q}^-(5,q^2)$
(see \cite[Theorem 8]{BKLP}).


\section{Partial quadrangles obtained from a hemisystem}\label{GQminushemi}

Recall that a \textit{hemisystem} $\mathcal{H}$ of a generalised quadrangle
$\mathcal{G}$ of order $(s,s^2)$, $s$ odd, is a set of $(s^3+1)(s+1)/2$ points
of $\mathcal{G}$ such that every line of $\mathcal{G}$ is incident with exactly
$(s+1)/2$ elements of $\mathcal{H}$. From $\mathcal{H}$, we construct a partial
quadrangle $\mathsf{PQ}(\mathcal{H})$ as follows:
\begin{center}
\begin{tabular}{l|l}
\textsc{points}& the points of $\mathcal{H}$\\
\textsc{lines}& the lines of $\mathcal{G}$.
\end{tabular}
\end{center}
The parameters are thus $((s-1)/2,s^2,(s-1)^2/2)$.  Since the complement of a
hemisystem is again a hemisystem, we may regard this construction as removing
``infinity'', where ``infinity'' is a hemisystem.

\begin{lemma}\label{glossaryHemi}
  Let $\mathcal{G}$ be a generalised quadrangle of order $(s,s^2)$, let
  $\mathcal{H}$ be a hemisystem of $\mathcal{G}$, and let $\I$ be an intriguing
  set of the partial quadrangle $\mathsf{PQ}(\mathcal{H})$ with intersection
  numbers $(h_1',h_2')$. Then we have the following information:

\begin{table}[H]
\begin{tabular}{l|c|c}
Case&Eigenvalue&Size\\
\hline
Negative intriguing set&$(-s^2+s-2)/2$&$h_2'(s+1)$\\
Positive intriguing set&$s-1$&$h_2'(s^3+1)/(s-1)^2$\\
Point set&$(s-1)(s^2+1)/2$&$(s+1)(s^3+1)/2$\\
\hline
\end{tabular}
\caption{Eigenvalues and sizes of intriguing sets of $ \mathsf{PQ}(\mathcal{H})$. }\label{table:glossaryHemi}
\end{table}
\end{lemma}

\begin{theorem}\label{atinfinity2}
  Let $\mathcal{G}$ be a generalised quadrangle of order $(s,s^2)$ and let
  $\mathcal{H}$ be a hemisystem of $\mathcal{G}$.  Let $\I$ be an intriguing set
  of $\mathcal{G}$ with parameters $(h_1,h_2)$. If $\I\setminus \mathcal{H}$ is
  an intriguing set of the partial quadrangle $\mathcal{G}\setminus\mathcal{H}$,
  then we have the following possibilities for the intersection numbers
  $(a_1,a_2)$ at infinity:
\begin{table}[H]
\begin{tabular}{l||c|c|c}
&$a_1-a_2$&$a_2$&$|\I\setminus\mathcal{H}|$\\
\hline
$-\rightarrow -$&$-(s^2+s)/2$&--&$(m(s^2+1)-a_2)(s+1)$\\
$-\rightarrow +$&$-(s^2+s)$&--&$(m(s^2+1)-a_2)(s^3+1)/(s-1)^2$\\
$+\rightarrow -$&$(s^2+s)/2$&$i/2$&$i(s+1)/2$\\
\hline
\end{tabular}
\caption{Possibilities for intersection numbers $(a_1,a_2)$.}\label{table:atinfinity2}
\end{table}
\end{theorem}

\begin{proof}
  The positive eigenvalues for $\mathcal{G}$ and
  $\mathcal{G}\setminus\mathcal{H}$ are both equal to $s-1$, however the
  negative eigenvalues differ: $-s^2-1$ for $\mathcal{G}$ and $(-s^2+s-2)/2$ for
  $\mathcal{G}\setminus \mathcal{H}$. Now we must have that $h_1-h_2$ and
  $(h_1-a_1)-(h_2-a_2)$ are eigenvalues for the respective geometries:

\begin{table}[H]
\begin{tabular}{ll|c|c||c|c|c}
Case&&$h_1-h_2$&$(h_1-a_1)-(h_2-a_2)$&$a_1-a_2$&$a_2$&$|\I\setminus\mathcal{H}|$\\
\hline
(i)&$-\rightarrow -$&$-s^2-1$&$(-s^2+s-2)/2$&$-(s^2+s)/2$&--&$(m(s^2+1)-a_2)(s+1)$\\
(ii)&$-\rightarrow +$&$-s^2-1$&$s-1$&$-(s^2+s)$&--&$(m(s^2+1)-a_2)(s^3+1)/(s-1)^2$\\
(iii)&$+\rightarrow -$&$s-1$&$(-s^2+s-2)/2$&$(s^2+s)/2$&$i/2$&$i(s+1)/2$\\
(iv)&$+\rightarrow +$&$s-1$&$s-1$&$0$&$i/2\frac{s^2-1}{s^2-s+1}$&$i(s+1)/2$\\
\hline
\end{tabular}
\caption{Details on the intersection numbers.}\label{table:casessofar}
\end{table}

\begin{enumerate}
\item[(i)] Suppose that $\I$ is an $m$-ovoid and $\I\setminus\mathcal{H}$ is
  negative intriguing. Then $|\I|=m(s^3+1)$ and hence
$$|\I\setminus\mathcal{H}|=\frac{m(s^2+1)-a_2}{(s-1)(s^2+1)-(-s^2+s-2)}(s+1)(s^3+1)
=(m(s^2+1)-a_2)(s+1).$$

\item[(ii)] Suppose that $\I$ is an $m$-ovoid and $\I\setminus\mathcal{H}$ is
  positive intriguing. Again we have $|\I|=m(s^3+1)$, but now we obtain
$$|\I\setminus\mathcal{H}|=\frac{m(s^2+1)-a_2}{(s-1)^2}(s^3+1).$$

\item[(iii)] Suppose that $\I$ is an $i$-tight set and $\I\setminus\mathcal{H}$
  is negative intriguing.  Then
$$|\I\setminus\mathcal{H}|=\frac{i-a_2}{(s-1)(s^2+1)-(-s^2+s-2)}(s^3+1)(s+1)
=(i-a_2)(s+1).$$ Since $\mathcal{H}$ is a hemisystem, we have that
$|\mathcal{H}\cap\I|=(s+1)i/2$.  So
$|\I\setminus\mathcal{H}|=|\I|-|\mathcal{H}\cap\I|=(s+1)i/2$. This gives
$(i-a_2)(s+1)=(s+1)i/2$ and hence $a_2=i/2$.

\item[(iv)] Suppose that $\I$ is an $i$-tight set and $\I\setminus\mathcal{H}$
  is positive intriguing.  Then $|\I|=i(s+1)$ and hence
$$|\I\setminus\mathcal{H}|=\frac{i-a_2}{(s-1)(s^2+1)-2(s-1)}(s^3+1)(s+1)
=\frac{i-a_2}{(s-1)^2}(s^3+1).$$ Since $\mathcal{H}$ is a hemisystem, we have
$|\mathcal{H}\cap\I|=(s+1)i/2$.  So
$|\I\setminus\mathcal{H}|=|\I|-|\mathcal{H}\cap\I|=(s+1)i/2$ and we obtain
$\frac{i-a_2}{(s-1)^2}(s^3+1)=(s+1)i/2$ and therefore
$a_1=a_2=i/2\frac{s^2-1}{s^2-s+1}$.  If we now compare with Lemma
\ref{glossaryHemi}, we arrive at the equation
$$2h_2'(s^3+1)=i(s-1)(s^2-1).$$
However, since we know $a_2$, and $h_2'=h_2-a_2$, we can substitute $h_2'$ in
the above equation and we obtain
$$2(1-\frac{s^2-1}{2(s^2-s+1)})(s^3+1)=(s-1)(s^2-1)$$
which implies that $2s=0$; a contradiction.
\end{enumerate}

\EndResult\end{proof}

\begin{corollary}\label{Q3qHemi}
  Let $\mathcal{G}$ be a generalised quadrangle of order $(s,s^2)$, $s$ odd, and
  let $\mathcal{H}$ be a hemisystem of $\mathcal{G}$. Both an $(s+1)$-tight set
  and an $(s^2+1)$-tight set of $\mathcal{G}$ never yield intriguing sets of
  $\mathsf{PQ}(\mathcal{H})$.  Furthermore, $\mathsf{Q}^+(3,q)$ and
  $\mathsf{Q}(4,q)$ embedded in $\mathsf{Q}^-(5,q)$, never yield intriguing sets
  of $\mathsf{PQ}(\mathcal{H})$.
\end{corollary}

\begin{proof}
  A $\mathsf{Q}^+(3,q)$ embedded in $\mathsf{Q}^-(5,q)$ is a $(q+1)$-tight set
  and a $\mathsf{Q}(4,q)$ section is a $(q^2+1)$-tight set.  Suppose
  $\mathcal{I}$ is an $(s+1)$-tight set of $\mathcal{G}$ such that
  $\mathcal{I}\setminus\mathcal{H}$ is a negative intriguing set of
  $\mathsf{PQ}(\mathcal{H})$ (the only case allowed by Theorem
  \ref{atinfinity2}).  Then we observe immediately a contradiction because
  $h_1-a_1$ is negative.  In the case that $\mathcal{I}$ is an $(s^2+1)$-tight
  set of $\mathcal{G}$, the parameter $a_1$ is equal to $s^2+(s+1)/2$. Now $a_1$
  is the number of points of $\mathcal{I}\cap \mathcal{H}$ which are collinear
  with an arbitrary point of $\mathcal{I}\setminus\mathcal{H}$. Therefore, $a_1$
  is divisible by $(s+1)/2$, which implies that $s+1$ divides
  $2s^2+s+1=(2s-1)(s+1)+2$; a contradiction.  So an $(s+1)$-tight set and an
  $(s^2+1)$-tight set of $\mathcal{G}$ never induce intriguing sets of
  $\mathsf{PQ}(\mathcal{H})$.
\EndResult\end{proof}

Examples exist for the first and third cases of Theorem \ref{atinfinity2} which
we demonstrate in what follows.  For the second case of Theorem
\ref{atinfinity2}, we do not have any examples when the generalised quadrangle
is $\mathsf{Q}^-(5,s)$. In this generalised quadrangle, an $m$-ovoid is a
hemisystem and we believe that only negative intriguing sets can arise in the
partial quadrangle.

\begin{conjecture}\label{hemi_negint2}\ \\
  Let $\mathcal{G}$ be a generalised quadrangle of order $(s,s^2)$ and let
  $\mathcal{H}$ be a hemisystem of $\mathcal{G}$.  Let $\I$ be another
  hemisystem of $\mathcal{G}$. Then $\I\setminus \mathcal{H}$ is a negative
  intriguing set of the partial quadrangle $\mathcal{G}\setminus\mathcal{H}$.
\end{conjecture}

The authors are not aware of a situation in which the above situation is
violated, and a proof of this fact would be a surprising result on the nature of
hemisystems.

\begin{lemma}[Nice Cone$\longrightarrow$ Negative Intriguing]\ \\
  Let $\mathcal{G}$ be a generalised quadrangle of order $(s,s^2)$, let
  $\mathcal{H}$ be a hemisystem of $\mathcal{G}$, and let $Z$ be a point of
  $\mathcal{H}$.  If the complement of $\mathcal{H}$ is intriguing at infinity
  for the partial quadrangle $\mathcal{G}\setminus Z^\perp$, then
  $Z^\perp\setminus \mathcal{H}$ is a negative intriguing set with parameters
  $((s-1)/2,(s^2+1)/2)$ of the partial
  quadrangle $\mathcal{G}\setminus\mathcal{H}$
\end{lemma}

\begin{proof}
  Let $X\in Z^\perp\setminus \mathcal{H}$. The points of $Z^\perp$ collinear
  with $X$ lie on the line $ZX$, and this line meets $\mathcal{H}'$ in $(s+1)/2$
  points. So $X$ is collinear with precisely $h_1'=(s-1)/2$ other points of
  $Z^\perp\setminus \mathcal{H}$. Now suppose that $X\notin Z^\perp\setminus
  \mathcal{H}$. So in particular $X$ is not in $Z^\perp$ and hence we can use
  the fact that $\mathcal{G}\setminus Z^\perp$ is a partial quadrangle.  Since
  $\mathcal{H}'$ is intriguing at infinity for the partial quadrangle
  $\mathcal{G}\setminus Z^\perp$, there exists a constant $a_2$ such that
  $|X^\perp\cap \mathcal{H}'\cap Z^\perp|=a_2$.  Now by Lemma \ref{hemi_negint},
  this value of $a_2$ is $(s^2+1)/2$.
\EndResult\end{proof}

If we have a partial spread $\mathcal{S}$ of a generalised quadrangle
$\mathcal{G}$, and a point $X$ not covered by any line of $\mathcal{S}$, then
$X$ is collinear with exactly one point of each member of $\mathcal{S}$. If half
of these points of collinearity are contained in a hemisystem, then we might
obtain an intriguing set of the associated partial quadrangle.

\begin{lemma}[Nice Partial Spread$\longrightarrow$ Negative Intriguing]\label{nicepartialspreads}\ \\
  Let $\mathcal{G}$ be a generalised quadrangle of order $(s,s^2)$ and let
  $\mathcal{H}$ be a hemisystem of $\mathcal{G}$.  Let $\I$ be a set of points
  covered by a partial spread of $c$ lines of $\mathcal{G}$ where $c$ is even.
  If $\I$ is intriguing at infinity, then for every point $X$ not in $\I$, half
  of the $c$ points of $\I$ collinear with $X$ are contained in $\mathcal{H}$,
  and $\I\setminus \mathcal{H}$ is a negative intriguing set of the partial
  quadrangle $\mathcal{G}\setminus\mathcal{H}$ with parameters
  $((c-s^2+s)/2-1,c/2)$.
\end{lemma}

\begin{proof}
It follows from Theorem \ref{atinfinity2}.
\EndResult\end{proof}

Such partial spreads as those described in Lemma \ref{nicepartialspreads} have
been found by computer for small $q$ in $\mathsf{Q}^-(5,q)$.  Therefore we have
examples in which a tight set of the generalised quadrangle induces negative
intriguing sets of the partial quadrangle arising from a hemisystem.

\begin{remark}
  A positive intriguing set with $h_2'=1$ has size $(s^3+1)/(s-1)^2$ which is
  only an integer when $s=2,3$.  Now let $\mathcal{O}$ be a maximal partial
  ovoid of a generalised quadrangle $\mathcal{G}$ of order $(s,s^2)$. If
  $\mathcal{O}$ is an intriguing set with parameters $(h_1',h_2')$ of
  $\mathsf{PQ}(\mathcal{H})$, then $h_1'=0$ and hence $\mathcal{O}$ is a
  negative intriguing set with $h_2'=(s^2-s+2)/2$.  In this case, we would have
  $|\mathcal{O}|=(s^2-s+2)(s+1)/2=(s^3+s+2)/2$, which happens to be the
  theoretical upper-bound for the size of a partial ovoid of $\mathsf{Q}^-(5,s)$
  given by De Beule, Klein, Metsch, and Storme \cite{DBKMS}.  In
  $\mathsf{Q}^-(5,3)$, there is a set of points of size $16$ that is a negative
  intriguing set of $\mathsf{PQ}(\mathcal{H})$ where $\mathcal{H}$ is Segre's
  hemisystem. It happens to be the unique maximal partial ovoid of
  $\mathsf{Q}^-(5,3)$ first discovered by Ebert and Hirschfeld
  \cite{EbertHirschfeld}.
\end{remark}

%
%

\section{Partial quadrangles that have a linear representation}\label{linear}

A \textit{$k$-cap} of a projective space $\mathsf{PG}(n,q)$ is a set of $k$
points with no three collinear.  Calderbank \cite{Calderbank82} proved using
number-theoretic arguments that if a partial quadrangle is a linear
representation then $q\ge 5$ or it is isomorphic to the linear representation of
one of the following:
\begin{enumerate}
\item An ovoid $\mathcal{O}$ of $\mathsf{PG}(3,q)$;
\item A Coxeter $11$-cap of $\mathsf{PG}(4,3)$;
\item A Hill $56$-cap of $\mathsf{PG}(5,3)$;
\item A $78$-cap of $\mathsf{PG}(5,4)$;
\item A $430$-cap of $\mathsf{PG}(6,4)$.
\end{enumerate}
Tzanakis and Wolfskill \cite{TzanakisWolfskill} then proved that if $q\ge 5$, we
must be in the first case. Since the only known ovoids are elliptic quadrics and
Suzuki-Tits ovoids, the examples in the first case are equivalent to the partial
quadrangles obtained from removing a point from a generalised quadrangle (of
order $(q,q^2)$). Hence we have just three known \textit{exceptional partial
  quadrangles} arising from (i) the Coxeter $11$-cap (yielding a
$\mathsf{PQ}(2,10,2)$), (ii) the Hill $56$-cap (yielding a
$\mathsf{PQ}(2,55,20)$) and (iii) the so-called Hill $78$-cap (yielding a
$\mathsf{PQ}(3,77,14)$). (It is still an open problem whether there exists a
$430$-cap of $\mathsf{PG}(6,4)$ or not.) For more details on these caps, we
refer the reader to Hill's paper \cite{Hill76}.

\begin{lemma}[Hyperplane $\longrightarrow$ Intriguing]\label{hyperplane}
  Let $\mathcal{K}$ be a cap of $\mathsf{PG}(n,q)$, and embed this projective
  space as a hyperplane $\pi_\infty$ of $\mathsf{PG}(n+1,q)$ so that the affine
  points and the lines meeting $\pi_\infty$ in a point of $\mathcal{K}$, form a
  partial quadrangle.  Let $\pi$ be the set of points in some hyperplane of
  $\mathsf{PG}(n+1,q)$ different from $\pi_\infty$ . Then $\pi$ is an intriguing
  set with parameters $$((q-1)|\pi\cap
  \mathcal{K}|,|\mathcal{K}\setminus\pi|).$$
\end{lemma}

\begin{proof}
  Let $P$ be a point of $\pi$. Then for every point $Q$ of $\pi\cap\mathcal{K}$,
  there are $q-1$ affine points on $QP$, other than $P$, which are collinear
  with $P$. Hence in total we have $(q-1)|\pi\cap \mathcal{K}|$ other points of
  the partial quadrangle collinear with $P$. Now suppose $P$ is not in
  $\pi$. Then clearly a point of $\mathcal{K}\cap \pi$ is not on a line
  connecting $P$ with a point of $\pi$. Since every line not in $\pi$ must meet
  $\pi$ in a point, it follows that the intersection number is
  $|\mathcal{K}\setminus\pi|$ in this case.
\EndResult\end{proof}

The example in the lemma above could either be a negative or positive intriguing
set depending on the intersection of the given hyperplane with the cap.

\begin{lemma}\label{glossaryexceptional}
  Let $\I$ be an intriguing set with intersection numbers $(h_1',h_2')$ of one
  of the three exceptional partial quadrangles. Then we have the following
  information:

\begin{table}[H]
\begin{tabular}{l|cc|cc|cc}
&Coxeter $11$-cap&&Hill $56$-cap&&Hill $78$-cap\\
Case&Eigenvalue&Size&Eigenvalue&Size&Eigenvalue&Size\\
\hline
Negative intriguing set&$-5$&$9h_2'$&$-23$&$(27/5)h_2'$&$-22$&$16h_2'$\\
Positive intriguing set&$4$&$(27/2)h_2'$&$4$&$(27/4)h_2'$&$10$&$(128/7)h_2'$\\
Point set&$22$&$243$&$112$&$729$&$234$&$4096$\\
\hline
\end{tabular}
\caption{Eigenvalues and sizes for intriguing sets of the exceptional partial quadrangles.}\label{table:glossaryexceptional}
\end{table}
\end{lemma}

By Lemma \ref{hyperplane}, the affine points in a hyperplane will have
associated eigenvalue $q|\pi\cap \mathcal{K}|-|\mathcal{K}|$ and so
$|\pi\cap\mathcal{K}|$ is $2$ or $5$ for the Coxeter $11$-cap, $11$ or $20$ for
the Hill $56$-cap, and $14$ or $22$ for the Hill $78$-cap.

\begin{lemma}[Nice Secundum $\longrightarrow$ Positive intriguing]\label{secunda}
  Let $\mathcal{K}$ be a cap of $\mathsf{PG}(n,q)$, and embed this projective
  space as a hyperplane at infinity of $\mathsf{PG}(n+1,q)$ so that the affine
  points and the lines meeting infinity in a point of $\mathcal{K}$, form a
  partial quadrangle.  Let $S$ be a secundum of $\mathsf{PG}(n+1,q)$ such that
  every hyperplane $\pi$ containing $S$ meets $\mathcal{K}$ in a constant number
  of points. Then the affine points of $S$ form an intriguing set with
  parameters $$((q-1)|S\cap \mathcal{K}|,|S\cap\pi|-|S\cap\mathcal{K}|).$$
\end{lemma}

\begin{proof}
  Let $X$ be a point of $S$, and let $C$ be a point of $\mathcal{K}$. If
  $C\notin S$, then there are no affine points of $S$ incident with $XC$, but if
  $C\in S$, then the affine points on the line $XC$ are all in $S$. So
  regardless of how $S$ meets the cap $\mathcal{K}$, it is clear that there are
  $1+(q-1)|\mathcal{K}\cap S|$ points collinear with $X$ in the associated
  partial quadrangle. So our first parameter is $(q-1)|\mathcal{K}\cap S|$. Now
  we look to the case that $X$ is not a point of $S$, and again, let $C$ be a
  point of $\mathcal{K}$.  Clearly $XC$ is not contained in $S$, but it may be
  disjoint from $S$ or meet $S$ in a point.  Let $\pi$ be the hyperplane
  $XS$. Now if $C\in S$, then $XC$ cannot meet $S$ in another point since
  otherwise $XC$ would be contained in $S$. If $C$ were not in $\pi$, then the
  unique point of intersection of $XC$ with $\pi$ would be $X$, and hence $XC$
  would not contain any points of $S$. So suppose $C\in \pi\setminus S$. Now
  $XC$ is a line of $\pi$, and $S$ is a hyperplane of $\pi$, thus $XC$ meets $S$
  in a point. Moreover, it is clear that this point of intersection is an affine
  point, and so the lines $XC$ which meet $S$ in a point are precisely those for
  which $C\in (\pi\cap \mathcal{K})\setminus(S\cap\mathcal{K})$. Hence, the
  affine points of $S$ form an intriguing set with second parameter equal to
  $|S\cap\pi|-|S\cap\mathcal{K}|$.
\EndResult\end{proof}

We remark that there are secunda of $\mathsf{PG}(5,3)$ which meet the Coxeter
$11$-cap in $3$ points, and hence every hyperplane containing such a secundum
must meet the Coxeter $11$-cap in $5$ points.  Similarly, there are secunda of
$\mathsf{PG}(6,3)$ which meet the Hill $56$-cap in $8$ points, and such that
every incident hyperplane meets this cap in $20$ points. Finally, we also have
secunda of $\mathsf{PG}(6,4)$ for the Hill $78$-cap which satisfy the hypotheses
of Lemma \ref{secunda}. Below we give some other examples which were found by
computer.


\subsection{Coxeter $11$-cap}

The permutation group induced on the Coxeter $11$-cap is $M_{11}$, and the full
stabiliser of the cap in $\mathsf{PGL}(6,3)$ is $3^5:(M_{11}\times 2)$.  We note
that this group is also the full automorphism group of the associated partial
quadrangle.  There were many negative intriguing sets found by computer, and we
report on those which were deemed interesting. There is a negative intriguing
set of size $45$ admitting $M_{10}$, and it is thus far, the only known negative
intriguing set of this size. Similarly, there are only two known negative
intriguing sets of size $54$, admitting groups of size $108$ and $864$
respectively. There is an intriguing set of size $81$ which is the complement of
the union of three hyperplanes (with stabiliser of size $648$). There are at
least two copies of $M_9:2$ in the automorphism group; one meets the normal
elementary abelian subgroup $3^5$ trivially, the other in a subgroup of order
$3^2$. These two groups give rise to intriguing sets of size $63$ and $108$, the
former is the complement of the disjoint union of negative intriguing sets of
size $45$.

There is a positive intriguing set of size $27$ which is the complement of the
union of $11$ hyperplanes, each meeting the cap in $5$ points. Its stabiliser is
$D_{18}\times S_3$. As noted above, there are solids of $\mathsf{PG}(5,3)$
meeting the cap in $3$ points, and hence we have positive intriguing sets of
size $27$.  All known examples arise from a sequence of unions and complements
of elements in the orbits of these two examples of size $27$.


\subsection{Hill $56$-cap}

The permutation group induced on the Hill $56$-cap is $\mathsf{PSL}(3,4).2$, and
the full stabiliser of the cap in $\mathsf{PGL}(7,3)$ is
$3^6:(2.\mathsf{PSL}(3,4).2)$.  We note that this group is also the full
automorphism group of the associated partial quadrangle.  The only known
negative intriguing set found so far is the set of affine points contained in a
hyperplane meeting the cap in 11 points.  As for positive intriguing sets, we
have hyperplanes on 20 cap points, solids on 4 cap points, and thousands of
other examples which are too numerous to list here. Most of these had
stabilisers of order $27$ or $54$.


\subsection{Hill $78$-cap}

The permutation group induced on the Hill $78$-cap is $(13:6)\times C_3$, and
the full stabiliser of the cap in $\mathsf{P\Gamma L}(8,4)$ is
$3^7:((C_{117}:C_3):2)$.  We note that this group is also the full automorphism
group of the associated partial quadrangle. Probably due to the fact that this
partial quadrangle has less symmetry than the other examples above, there were
many intriguing sets found, and none believed to be particularly interesting to
the authors.  The smallest negative intriguing set found had size $512$ (so with
parameters $(10,32)$), and the smallest positive intriguing set had size $128$
(parameters $(17,7)$) and hence attains the minimum size. Most of the intriguing
sets found had their full stabiliser acting regularly on them.

\section{Concluding remarks}

We introduced the definition of an intriguing set via strongly regular graphs,
and although much of the interest so far has been on intriguing sets of
generalised quadrangles and partial quadrangles, it may perhaps also be
interesting to investigate the intriguing sets of other particular families of
strongly regular graphs.

\section*{Acknowledgements}
This work was supported by the GOA-grant ``Incidence Geometry'' at Ghent
University. The first author acknowledges the support of a Marie Curie Incoming
International Fellowship within the 6th European Community Framework Programme
(contract number: MIIF1-CT-2006-040360), and the third author was supported by a
travel fellowship from the School of Sciences and Technology -- University of
Naples ``Federico II''.

\end{document}